\documentclass[letterpaper]{article} 
\usepackage[]{aaai23}  
\usepackage{times}  
\usepackage{helvet}  
\usepackage{courier}  
\usepackage[hyphens]{url}  
\usepackage{graphicx} 
\def\UrlFont{\rm}  
\usepackage{natbib}  
\usepackage{caption} 
\usepackage{algorithm}
\usepackage{algorithmic}

\usepackage{newfloat}
\usepackage{listings}
\DeclareCaptionStyle{ruled}{labelfont=normalfont,labelsep=colon,strut=off} 
\floatstyle{ruled}
\newfloat{listing}{tb}{lst}{}
\floatname{listing}{Listing}
\usepackage[utf8]{inputenc} 
\usepackage[T1]{fontenc}    
\usepackage{url}            
\usepackage{booktabs}       
\usepackage{amsfonts}       
\usepackage{nicefrac}       
\usepackage{microtype}      
\usepackage{xcolor}         

\usepackage{graphicx}
\usepackage{amssymb}
\usepackage{amsmath}
\usepackage{amsthm}
\usepackage{ragged2e} 
\usepackage{stackengine}
\usepackage{bbm}
\usepackage{enumitem} 
\usepackage{bbm}

\newcommand{\rvline}{\hspace*{-\arraycolsep}\vline\hspace*{-\arraycolsep}}

\usepackage[english]{babel}
\newtheorem{theorem}{Theorem}

\newtheorem{definition}{Definition}
\newtheorem{proposition}{Proposition}

\newtheorem{lemma}{Lemma}
\newtheorem{example}{Example}

\usepackage{mathtools}				
\DeclareMathOperator{\rk}{rank}	
\DeclareMathOperator{\tr}{tr}	    
\DeclareMathOperator{\vecc}{vec}	    
\DeclareMathOperator{\mat}{mat}	    
\DeclareMathOperator{\st}{s.t.}

\newcommand{\RR}{\mathbb R}
\DeclareMathOperator{\RIP}{RIP}

\DeclarePairedDelimiter{\norm}{\lVert}{\rVert}

\title{Semidefinite Programming versus Burer-Monteiro Factorization for Matrix Sensing}

\author{
    Baturalp Yalcin,\equalcontrib\textsuperscript{\rm 1}
    Ziye Ma,\equalcontrib\textsuperscript{\rm 2}
    Javad Lavaei, \textsuperscript{\rm 1}
    Somayeh Sojoudi\textsuperscript{\rm 2}
}
\affiliations {
    \textsuperscript{\rm 1}UC Berkeley, Department of Industrial Engineering and Operations Research,\\
    \textsuperscript{\rm 2}UC Berkeley, Department of Electrical Engineering and Computer Science\\
    \{baturalp\textunderscore yalcin, ziyema, lavaei, sojoudi\}@berkeley.edu
}

\nocopyright

\begin{document}

\maketitle

\begin{abstract}
    Many fundamental low-rank optimization problems, such as matrix completion, phase synchronization/retrieval, power system state estimation, and robust PCA, can be formulated as the matrix sensing problem. Two main approaches for solving matrix sensing are based on semidefinite programming (SDP) and Burer-Monteiro (B-M) factorization. The SDP method suffers from high computational and space complexities, whereas the B-M method may return a spurious solution due to the non-convexity of the problem. The existing theoretical guarantees for the success of these methods have led to similar conservative conditions, which may wrongly imply that these methods have comparable performances. In this paper, we shed light on some major differences between these two methods. First, we present a class of {\it structured} matrix completion problems for which the B-M methods fail with an overwhelming probability, while the SDP method works correctly. Second, we identify a class of {\it highly sparse} matrix completion problems for which the B-M method works and the SDP method fails. Third, we prove that although the B-M method exhibits the same performance independent of the rank of the unknown solution, the success of the SDP method is correlated to the rank of the solution and improves as the rank increases. Unlike the existing literature that has mainly focused on those instances of matrix sensing for which both SDP and B-M work, this paper offers the first result on the unique merit of each method over the alternative approach. 
    
\end{abstract}

 \section{Introduction}

 Low-rank matrix recovery problems have ubiquitous applications in machine learning and data analytics, including collaborative filtering \citep{koren2009matrix}, phase retrieval \citep{candes2015phase,singer2011angular,boumal2016nonconvex,shechtman2015phase}, motion detection \citep{fattahi2020exact}, and power system state estimation \citep{jin2020boundary,zhang2017conic,jin2019towards}. This problem is formally defined as follows: Given a measurement operator $\mathcal{A}(\cdot): \mathbb{R}^{m \times n} \mapsto \mathbb{R}^{d}$ returning a $d$-dimensional measurement vector $\mathcal{A}(\mathbf{M}^*)$ from a low-rank ground truth matrix $\mathbf{M}^* \in \mathbb{R}^{m \times n}$ with rank $r$, the goal is to obtain a matrix with rank less than equal to $r$ that conforms with the measurements, preferably the ground truth matrix $\mathbf{M^*}$. This problem can be stated as the feasibility problem
\begin{align} \label{eqn:gen}
    \mathrm{find}&\quad {\bf M}\in\mathbb{R}^{m\times n} \\
    \mathrm{s.t.}& \quad \mathcal{A}(\mathbf{M}) = \mathcal{A}(\mathbf{M}^*) \notag\\
    & \quad \rk({\bf M}) \leq r.
    \notag
\end{align}
While the measurement operator $\mathcal{A}$ can be nonlinear as in the case of one-bit matrix sensing \citep{davenport20141} and phase retrieval \citep{shechtman2015phase}, matrix sensing and matrix completion that are widely studied have linear measurement operators \citep{candes2009exact,recht2010guaranteed}. We focus on the matrix sensing and matrix completion problems throughout this paper. Despite the linearity of $\mathcal{A}$, there are two types of problems depending on the structure of the ground truth matrix $\mathbf{M}^*$. The first type, symmetric problem, consists of a low-rank positive semidefinite ground truth matrix $\mathbf{M}^* \in \mathbb{R}^{n \times n}$, whereas the second type, asymmetric problem, consists of a ground truth matrix $\mathbf{M}^* \in \mathbb{R}^{m \times n}$ that is possibly sign indefinite and non-square. Since each asymmetric problem can be converted to an equivalent symmetric problem \citep{zhang2021general}, we study only the symmetric problem in this paper.

The matrix sensing and completion problems have linear measurements; hence, the first constraint in problem \eqref{eqn:gen} is linear. Therefore, the only nonconvexity of the problem arises from the nonconvex rank constraint. Earlier works on these problems focused on their convex relaxations by penalizing high-rank solutions \citep{candes2009exact,recht2010guaranteed, candes2010power}. They utilized the nuclear norm of a matrix as the convex surrogate of the rank function. This led to semidefinite programming (SDP) relaxations, which solve the original non-convex problems exactly with high probability based on some assumptions on the linear measurement operator and the ground truth matrix, such as the Restricted Isometry Property (RIP) and incoherence conditions. High computational time and storage requirements of the SDP algorithms incentivized the implementation of the B-M factorization approach \citep{burer2003nonlinear}. This approach factorizes the symmetric matrix variable $\mathbf{M} \in \mathbb{R}^{n \times n} $ as $\mathbf{M} = \mathbf{XX}^T$ for some matrix $\mathbf{X} \in \mathbb{R}^{n \times r}$, which obviates imposing the positive semidefiniteness and rank constraints. Although the dimension of the decision variable reduces dramatically when $r$ is small, the problem is still nonconvex since its objective function is nonconvex in terms of the factorized $\mathbf{X}$.

\subsection{Problem Formulation}
Formally, the SDP formulation of the matrix sensing problem uses the nuclear norm of the variable, $\| \mathbf{M} \|_*$, to serve as a surrogate of the rank, and replaces the rank constraint in \eqref{eqn:gen} with an objective to minimize $\| \mathbf{M} \|_*$. Due to the symmetricity and positive semidefiniteness of the variable, the nuclear norm is equivalent to the trace of the matrix variable $\mathbf{M}$. Hence, the SDP formulation can be written as
\begin{align} \label{eqn:ms-mc-sdp}
    \underset{{\bf M} \in \mathbb{R}^{n \times n}}{\mathrm{min}} &\quad \mathrm{tr}({\bf M}) \qquad
    \mathrm{s.t.} \ \ \mathcal{A}(\mathbf{M}) = b, \  \mathbf{M} \succeq 0, 
\end{align}
where $b = \mathcal{A}(\mathbf{M}^*) = [ \langle \mathbf{A}_1, \mathbf{M}^* \rangle, \dots, \langle \mathbf{A}_d, \mathbf{M}^* \rangle]^T$ is given and $\{ \mathbf{A}_i \}_{i=1}^d \in \mathbb{R}^{n \times n}$ are called sensing matrices. Moreover, the matrix completion problem is a special case of the matrix sensing problem with each sensing matrix measuring only one entry of $\mathbf{M}^*$. We can represent the measurement operator $\mathcal{A}$ as $\mathcal{A}_\Omega : \mathbb{R}^{n \times n} \mapsto \mathbb{R}^{n \times n}$ for this special case, which is defined as follows: 
\[ \mathcal{A}_\Omega({\bf M})_{ij} := \begin{cases} {\bf M}_{ij} & \text{if } (i,j) \in \Omega\\ 0 & \text{otherwise}, \end{cases} \]
where $\Omega$ is the set of indices of observed entries. We denote the measurement operator as ${\bf M}_\Omega := \mathcal{A}_{\Omega}({\bf M})$ for simplicity.
Besides the SDP formulation, the B-M factorization formulation of the matrix sensing (MS) and matrix completion (MC) problems can be stated as
\begin{subequations}
    \begin{align} \label{eqn:ms-bm}
    & \text{(MS)} \quad  \min_{{\bf X}\in\mathbb{R}^{n\times r}}~  g\left[ \mathcal{A}({\bf XX}^T) - b \right],\\ 
   \label{eqn:mc-bm} & \text{(MC)} \quad  \min_{{\bf X}\in\mathbb{R}^{n\times r}}~ g\left[ (\mathbf{XX}^T  - {\bf M}^*)_\Omega \right],
\end{align}
\end{subequations}
where $g(\cdot): \mathbb{R}^d \mapsto \mathbb{R}$ is some twice continuously differentiable function such that $\mathbf{0}_{n \times n}$ is its unique minimizer and the Hessian of $g(\cdot)$ is positive definite at $\mathbf{0}_{n \times n}$. These assumptions are satisfied by the common loss functions considered in the literature. The main objective of this paper is to compare the SDP and B-M methods for the MC and MS problems. 

\subsection{Background and Related Work}
It is widely known that the SDP formulation \eqref{eqn:ms-mc-sdp} can be used to solve the matrix sensing problem if the sensing matrices are sampled independently from a sub-Gaussian distribution and the number of measurements $d$ is large enough \citep{recht2010guaranteed,recht2008necessary}. This is also a sufficient condition for the sensing matrices to satisfy the RIP condition with high probability, which is defined below:
\begin{definition}[RIP] \citep{candes2009exact} \label{def:rip}
    The linear map $\mathcal{A}: \mathbb{R}^{n \times n} \mapsto \mathbb{R}^{m}$ is said to satisfy $\delta_{p}$-RIP if there is a constant $\delta_{p} \in [0,1)$ such that
    \[ (1-\delta_{p})\|\mathbf{M} \|_F^2 \leq \| \mathcal{A} (\mathbf{M}) \|^2 \leq (1 + \delta_{p}) \|\mathbf{M}\|_F^2 \]
    holds for all matrices $\mathbf{M} \in \mathbb{R}^{n \times n}$ satisfying $\rk(\mathbf{M}) \leq p$.
\end{definition}
The RIP constant $\delta_{p}$ represents how similar the linear operator $\mathcal{A}$ is to an isometry, and various upper bounds on $\delta_p$ have been proposed to serve as sufficient conditions for the exact recovery (meaning that one can recover the ground truth $\bf M^*$ by solving the SDP problem). A few notable ones include $\delta_{4r} < \sqrt{2}-1$ in \citep{candes2010tight}, $\delta_{5r} < 0.607, \delta_{3r} < 0.472$ in \citep{mohan2010new}, and $\delta_{2r} < 1/2, \delta_{3r} < 1/3$ in \citep{cai2013sharp}. On the other hand, when the sensing matrices are not sampled independently from a sub-Gaussian distribution or when the RIP condition is not met, the SDP formulation may still recover the ground truth matrix with a high probability. This is the case for MC problems for which RIP fails to hold while SDP works as long as entries of observation follow an independent Bernoulli model \citep{candes2009exact,candes2010power}.

However, recent works have shown that if we use the B-M method instead of the SDP approach, we can still recover the ground truth matrix via first-order methods under similar RIP or coherence assumptions in both the matrix sensing and matrix completion cases \citep{ge2017no,bhojanapalli2016global,PKCS2017,ZWYG2018,ZLTW2018,zhang2019sharp,zhang2020many,bi2020global,HLB2020,zhu2021global,Zhang2021-p,ZBL2021,ma2021sharp,ma2022noisy}. Namely, the state-of-the-art result states that as long as $\delta_{\tilde r+r} < 1/2$ for the matrix sensing problem, there exists no spurious{\footnote{A local minimum is called spurious if it is not a global minimum. } local minima for an over-parametrized B-M formulation and the gradient descent algorithm can recover $\mathbf{M}^*$ exactly \citep{Zhang2021-p}. Here, $\tilde r \geq r$ is the search rank that we choose manually in the B-M formulation. If we know the value of $r$, we can set $\tilde r$ to $r$, making the B-M approach enjoy the same RIP guarantee as the SDP approach. Since the B-M approach enjoys far better scalability, it has become an increasingly popular tool for solving the matrix sensing problem.

Nevertheless, the B-M approach cannot be routinely used without careful consideration since it could fail on easy (from an information-theoretic perspective) instances of the problem as demonstrated in \citep{yalccin2022factorization}, especially in cases when the RIP condition is not satisfied. 

Thus, it is important to compare and contrast both the SDP and B-M approaches to discover which method is superior to the other one. This comparison is timely since specialized sparse SDP algorithms have become more efficient in recent years, making the SDP method more practical than before \citep{zhang2021sparse, udell2017, yurtsever21}. In this paper, we show that the SDP approach is more powerful than the B-M method as far as the RIP measure is concerned. We also discover that the B-M method is able to solve certain instances for which the SDP approach fails. This means that none of these techniques is universally better than the other one and the best technique should be chosen based on the nature of the problem. This work provides the first step towards understanding the trade-off between a well-known convex relaxation and first-order descent algorithms applied to the B-M factorization formulation.

\subsection{Our Contributions}

We provide a comparative analysis between the SDP approach and the B-M method. We first present the advantages of the SDP approach over the B-M method:
\begin{enumerate}
    \item First, we focus on an important class of MC problems recently studied in \cite{yalccin2022factorization}. That paper has shown that even though this class has low information-theoretic complexity, the B-M method would utterly fail and the probability of success via first-order methods is almost zero. We prove that the SDP method successfully solves this class and, therefore, SDP may not suffer from the unusual behavior of B-M with regard to easy instances of MC. This also implies that the information-theoretic and optimization complexities are expected to be more aligned for SDP than B-M.
    
    \item We then investigate a class of MS problems found in the recent paper \cite{zhang2022unified}. Each MS instance belonging to this class satisfies $\delta_2$-RIP with $r=1$ for some $\delta > 1/2$ such that the B-M formulation leads to $\mathcal{O}((1-\delta)^{-1})$ spurious solutions and this number goes to infinity as $\delta$ approaches $1$. We show that although each instance is extremely non-convex based on the number of spurious solutions, the SDP method successfully solves all of the problems in this class. This implies that, unlike the B-M method, the success of the SDP approach is not directly correlated to the presence of many spurious solutions.    

    \item The recent paper \cite{zhang2021general} has shown that the sharpest RIP bound for the success  of the B-M method on the MS problem is $1/2$ and this is independent of the rank $r$. This is an undesirable result since high-rank problems have lower information-theoretic complexity than low-rank problems. We derive a sufficient RIP bound for the SDP method and show that it can increase from $1/2$ to $1$ as the rank $r$ becomes larger. This implies that the SDP approach does not suffer from a major shortcoming of the B-M method.
\end{enumerate}
Despite the above advantages, we show that the SDP approach is not universally better than the B-M method. To prove this, we identify a class of MC problems with $\mathcal{O}(n)$ observations in the rank-$1$ case for which B-M works while SDP fails. It is clear from these comparisons that although the B-M approach is known to be more powerful due to its scalability property, the SDP approach enjoys some unique merits and deserves to be revisited, especially in light of the advancements of fast SDP solvers \cite{zhang2021sparse, udell2017, yurtsever21}

\section{Notations}

The symbol $[n]$ represents the set of integers from $1$ to $n$. We use lower-case bold letters, namely $\mathbf{x}$, to represent vectors and capital bold letters, namely $\mathbf{X}$, to represent matrices. $\mathbf{I}_n$ refers to the identity matrix of size $n \times n$ and $\mathbf{0}_{n \times n}$ refers to the $n \times n $ dimensional matrix with zero entries.  $\norm{\mathbf{x}}$ denotes the Euclidean norm of the vector $\mathbf{x}$, $\|\mathbf{X}\|$ and $\|\mathbf{X}\|_F$ are the $2$-norm and the Frobenius norm of the matrix $\mathbf{X}$, respectively. For every vector $\bf x$, $[\bf x]_i$ denotes the $i$-th entry and $[\bf x]_{i:j}$ denotes the subvector of entries from index $i$ to index $j$ for $i < j$. Similarly, for every matrix $\bf X$, $[\bf X]_{i:j, k:l}$ denotes the submatrix with rows between $i$ and $j$ and columns between $k$ and $l$ with $i<j$ and $k<l$. Let $\langle {\bf A}, {\bf B} \rangle = \tr({\bf A}^T{\bf B})$ be the inner product between matrices. The Kronecker product between $\bf A$ and $\bf B$ is denoted as $\bf A \otimes \bf B$. For a matrix $\bf X$, $\vecc(\bf X)$ is the usual vectorization operation by stacking the columns of the matrix $\bf X$ into a vector. For a vector $\mathbf{x} \in \RR^{n^2}$, $\mat(\bf x)$ converts $\bf x$ to a square matrix and $\mat_S(\bf x)$ converts $ \bf x$ to a symmetric matrix, i.e., $\mat( \bf x)= \bf X$ and $\mat_S(\bf x) = ( \bf X + \bf X^T)/2$, where $\bf X \in \RR^{n \times n}$ is the unique matrix satisfying $ \bf x=\vecc(\mathbf{X})$.  The notations ${\bf X}\succeq 0$ and ${\bf X}\succ 0$ mean that the matrix ${\bf X}$ is positive semidefinite (PSD) and positive definite, respectively. The set of $n \times n$ PSD matrices is denoted as $\mathbb{S}_+^n$.
For a function $f: \mathbb{R}^{m \times n} \mapsto \mathbb{R}$, we denote the gradient and the Hessian as $\nabla f(\cdot)$ and $\nabla^2f(\cdot)$, respectively. The Hessian is a four-dimensional tensor with $[\nabla^2 f(\mathbf{X})]_{i,j,k,l} = \frac{\partial^2 f(\mathbf{X})}{\partial \mathbf{X}_{i,j} \partial \mathbf{X}_{k,l}}$ for all $i,j \in [m]$ and $k,l \in [n]$.  We use $\lceil\cdot\rceil$ and $\lfloor\cdot\rfloor$ to denote the ceiling and floor operators, respectively. The cardinality of a set $\mathcal{S}$ is shown as $|\mathcal{S}|$.

\section{Advantages of the SDP Approach}


\subsection{B-M Fails While SDP Succeeds}
In this section, we focus on a class of MC instances that was first proposed in \citep{yalccin2022factorization} for which the B-M factorization fails. We focus on the matrix completion problem since it is the most common special case of the matrix sensing problem that does not satisfy the RIP condition. We will prove that while the B-M approach fails to recover $\mathbf{M}^*$, the SDP approach can provably find $\mathbf{M}^*$.

We will first give an introduction to this class of MC instances. Consider a rank-$r$ ground truth matrix $\mathbf{M}^* \in \mathbb{S}_{+}^n$ with $r \geq 1$ and $n \geq 2r$. Let $m := n/r$ and assume without the loss of generality that $n$ is divisible by $r$.  We decompose the ground truth matrix into blocks of dimension $r \times r$; thus, $\mathbf{M}^*$ is an $m \times m$ block matrix whose block element at the position $(i,j)$ is denoted as $\mathbf{M}^*_{i,j}$ for $i,j \in [m]$. We require some graph-theoretic notions before introducing the underlying class of MC instances.
\begin{definition}[Induced Measurement Set] \label{def:block-sparsity}
    Let $\mathcal{G} = (\mathcal{G}_1,\mathcal{G}_2)= (\mathcal{V}, \mathcal{E}_1, \mathcal{E}_2)$ be a pair of undirected graphs with the node set $\mathcal{V} = [m]$ and the disjoint edge sets $\mathcal{E}_1,\mathcal{E}_2  \subset [m]\times[m]$, respectively. The induced measurement set $\Omega(\mathcal{G})$ is defined as follows: if $(i,j) \in \mathcal{E}_1$, then the entire block ${\bf M}^*_{i,j}$ is observed; if $(i,j) \in \mathcal{E}_2$, then all nondiagonal entries of the block ${\bf M}^*_{i,j}$ are observed; otherwise, none of the entries of the block is observed. The graph $\mathcal{G}$ is referred to as the block sparsity graph.
\end{definition}
We represent the general problem \eqref{eqn:gen} with the linear measurement operator $\mathcal{A}$ and rank-$r$ ground truth matrix $\mathbf{M}^* \in \mathbb{R}^{n \times  n}$ as $\mathcal{P}_{\mathbf{M}^*, \mathcal{A}, n, r}$. If this is a matrix completion problem with the measurement set $\Omega$, then this special case of the same problem is denoted as $\mathcal{P}_{\mathbf{M}^*, \Omega, n, r}$.
Based on Definition \ref{def:block-sparsity} and this notation, a low-complexity class of MC instances will be introduced below. These instances have a low complexity because graph-theoretical algorithms can solve them in polynomial time in terms of $n$ and $r$.
\begin{definition}[Low-complexity class of MC instances] \label{def: low-comp}
    Define $\mathcal{L}(\mathcal{G}, n, r)$ to be the class of low-complexity MC instances $\mathcal{P}_{{\bf M}^*, \Omega, n, r}$ with the following properties:
    \setlist{nolistsep}
    \begin{enumerate}[label= \bf \roman*), noitemsep]
        \item The ground truth $\mathbf{M^*} \in \mathbb{S}^{n}_+$ is rank-r.
        \item The matrix $\mathbf{M}^*_{i,j}\in\mathbb{R}^{r\times r}$ is rank-$r$ for all $i,j\in[m]$.
        \item The measurement set $\Omega=\Omega(\mathcal{G})$ is induced by $\mathcal{G}=(\mathcal{G}_1,\mathcal{G}_2)$, where $\mathcal{G}_1$ is connected, non-bipartite, and its vertices have self-loops.
    \end{enumerate}
\end{definition}
The next theorem borrowed from \citep{yalccin2022factorization} illustrates the failure of the B-M factorization method.
\begin{theorem} \label{thm: old-paper}
    Consider a maximal independent set $\mathcal{S}(\mathcal{G}_1)$ of $\mathcal{G}_1$ such that the induced subgraph by vertices in $\mathcal{S}$, $\mathcal{G}_2[\mathcal{S}]$, is connected. There exists an instance in $\mathcal{L}(\mathcal{G}, n, r)$ for which the problem \eqref{eqn:mc-bm} has at least $2^{r|\mathcal{S}(\mathcal{G}_1)|}-2^r$ spurious local minima. In addition, the randomly initialized gradient descent algorithm converges to a global minimum with probability at most $\mathcal{O}(2^{-r|\mathcal{S}(\mathcal{G}_1)|})$, while there is a graph-theoretical algorithm that can solve the problem in $\mathcal{O}(n^2/r^2 + nr^2)$ time.
\end{theorem}
%

The proof of Theorem \ref{thm: old-paper} utilizes the Implicit Function Theorem (IFT). Specifically, the work \cite{yalccin2022factorization} has generated ground truth matrices $\mathbf{M}^*$ for which the B-M method has $2^{r|\mathcal{S}(\mathcal{G}_1)|}$ global solutions and only $2^r$ of them correspond to the correct completion of the $\mathbf{M}^*$. A generic small perturbation of the problem results in a new instance of an MC problem that belongs to the low-complexity class of MC instances. The conditions on $\mathcal{G}_1$ guarantee that the perturbed problem belongs to the low-complexity class, while the conditions on $\mathcal{G}_2$ guarantee that the Hessian of the objective function of the unperturbed problem is positive definite at the global solutions. Since the instances in the low-complexity class are well defined, the new perturbed problem has a unique completion with $2^r$ possible global solutions for the B-M method. On the other hand, the other stationary points that correspond to global solutions of the unperturbed problem must be spurious local minima of the new instance. This is concluded by using the IFT. The perturbation that yields a new instance in the low-complexity class of the MC problem is achieved by perturbing the ground truth matrix $\mathbf{M}^* = \mathbf{X}^*(\mathbf{X}^*)^T$ by a small and generic perturbation $\epsilon \in \mathbb{R}^{n \times r}$. The new ground truth matrix is $\mathbf{M}^*(\epsilon) = \mathbf{X}^*(\epsilon)(\mathbf{X}^*(\epsilon))^T$, where $\mathbf{X}^*(\epsilon)_i = \mathbf{X}^*_i + \epsilon_i$ if $i \in \mathcal{S}(\mathcal{G}_1)$ and $\mathbf{X}^*(\epsilon)_i = \epsilon_i$ otherwise and $\rk(\mathbf{X}^*_i) = \rk(\mathbf{X}^*_i + \epsilon_i) = r, \forall i \in [m]$. A generic perturbation $\epsilon$ does not belong to a measure zero set in $\mathbb{R}^{n \times r}$.

It is desirable to study how the SDP method performs on this low-complexity class of MC instances. We will present the result for a larger class of problems that contains all instances discussed in Theorem \ref{thm: old-paper}.
\begin{theorem} \label{thm: sdp-solves-mc}
    Given $\mathcal{G} = (\mathcal{G}_1, \mathcal{G}_2) = (\mathcal{V}, \mathcal{E}_1, \mathcal{E}_2)$, consider any maximal independent set $\mathcal{S}(\mathcal{G}_1)$. Consider also $\mathbf{M}^*(\epsilon) = \mathbf{X}^*(\epsilon)(\mathbf{X}^*(\epsilon))^T$ for any arbitrary $\epsilon \in \mathbb{R}^{n \times r}$, where $\mathbf{X}^*(\epsilon)_i = \mathbf{X}^*_i + \epsilon_i$ if $i \in \mathcal{S}(\mathcal{G}_1)$ and $\mathbf{X}^*(\epsilon)_i = \epsilon_i$ otherwise and $\rk(\mathbf{X}^*_i) = \rk(\mathbf{X}^*_i + \epsilon_i) = r, \forall i \in [m]$. The SDP formulation \eqref{eqn:ms-mc-sdp} with the observation set $\Omega$ induced by $\mathcal{G}_1$ uniquely recovers the ground truth matrix $\mathbf{M}^*(\epsilon)$.
\end{theorem}
Note that we do not require $\epsilon$ to be small or have access to partial observations of blocks induced by edges in $\mathcal{G}_2$. Hence, Theorem \ref{thm: sdp-solves-mc} shows that SDP solves all MC instances introduced in Theorem 1 and beyond. As a result of Theorem \ref{thm: sdp-solves-mc}, the SDP approach is a viable choice for those MC instances for which the preferable and faster B-M factorization method fails to recover the ground truth matrix. Similar to perturbing the ground truth matrix, one can perturb the linear measurement operator of the matrix completion problem $\mathcal{A}_\Omega$ as 
\begin{equation} \label{eqn: operator-ms}
    \mathcal{A}_{\Omega(\epsilon)}({\bf M})_{ij} := \begin{cases} {\bf M}_{ij}, & \text{if } (i,j) \in \Omega\\ \epsilon \mathbf{M}_{ij}, & \text{otherwise} \end{cases},
\end{equation}
where $\epsilon > 0$ is a sufficiently small real number \citep{zhang2022unified}. Note that $\mathcal{A}_{\Omega(\epsilon)}$ satisfies the RIP condition with $\delta = (1-\epsilon)/(1+\epsilon)$. 
\begin{theorem} \label{thm: bm-fails-ms}
    Suppose that $g$ is the squared loss function, i.e $g(\bf x) = \| \bf x \|^2$. Consider the measurement set $\Omega$ defined in Theorem \ref{thm: old-paper}. For every sufficiently small $\epsilon > 0$, there exists a low-complexity instance of the MS problem $\mathcal{P}_{\mathbf{M}^*, \mathcal{A}_{\Omega}(\epsilon), n, r}$ with $\mathcal{O}({2^{r|\mathcal{S}(\mathcal{G}_1) |}})$ spurious local minima.
\end{theorem}
The proof of the above theorem is similar to the proof of Theorem \ref{thm: old-paper} because the conditions for unperturbed problems are the same and a different small perturbation to the problem yields a similar number of spurious solutions. Hence, the proof is omitted. The above theorem states that there are not only MC instances but also MS instances that suffer from this undesirable behavior of the B-M factorization approach. The ground truth matrix $\mathbf{M}^*$ is generated as in Theorem \ref{thm: old-paper} to have $2^{r|\mathcal{S}(\mathcal{G}_1)|}$ global solutions for the unperturbed problem. Furthermore, the number of spurious solutions for this scheme can be quantified as $\mathcal{O}((1-\delta)^{-1})$ for $\delta \geq  1/2$ in the rank-$1$ case \citep{zhang2022unified}.  Nevertheless, the SDP formulation approach trivially solves all these undesirable MS instances. This is due to the fact the perturbed measurement operator $\mathcal{A}_{\Omega(\epsilon)}$ corresponds to observing all the entries. Hence, the feasible set only contains the ground truth matrix $\mathbf{M}^*$.
\begin{proposition} \label{prop:sdp-solves-ms}
    Given a measurement set $\Omega$, the SDP formulation \eqref{eqn:ms-mc-sdp} uniquely recovers the rank-r ground truth matrix $\mathbf{M}^*$ for the MS instance $\mathcal{P}_{\mathbf{M}^*, \mathcal{A}_{\Omega}(\epsilon), n, r}$, where $\mathcal{A}_{\Omega}(\epsilon)$ is defined in \eqref{eqn: operator-ms} and $\epsilon$ is an arbitrary nonzero number. 
\end{proposition}
Hence, the SDP approach successfully solves all the instances in Theorem \ref{thm: bm-fails-ms} for which the RIP constant exists (while greater than $1/2$), unlike the B-M method. Consequently, SDP could be the preferred method when the sufficient conditions on RIP for exact recovery by the B-M factorization are not met. In the next part, we will provide sharper sufficiency bounds for the SDP approach, which further corroborates its strength.

\subsection{Sharper RIP bound for SDP}
Since the SDP method is more powerful than the B-M factorization for certain classes of MC and MS problems as shown in the previous section and since specialized SDP algorithms can solve large-scale MC and MS problems, it is useful to further study the SDP method through the lens of the well-known RIP notion. We will derive a strong lower bound $\delta_{lb}$ on the RIP constant $\delta$ to guarantee convergence to the ground truth solution by using a proof technique called the inexistence of incorrect solution \cite{zhang2019sharp}. We aim to find a linear measurement operator $\mathcal{A}$ with the smallest RIP constant such that the SDP formulation converges to a wrong solution. To do so, we need to solve the optimization problem
\begin{equation} \label{eqn:sdp-inex}
    \begin{aligned}
    \min_{\delta,\mathcal{A}} \quad & \delta \\
    \text{s.t.} \quad & \mathcal{A}(\mathbf{M}) = \mathcal{A}(\mathbf{M}^*) \\
    & \tr(\mathbf{M}) \leq \tr(\mathbf{M}^*) \\
    & \text{$\mathcal{A}$ satisfies the $\delta-{2r}$-$\RIP$ property},
\end{aligned}
\end{equation}
where $\mathbf{M} \not = \mathbf{M}^*$. The condition $tr(\bf M)\leq tr(\bf M^*)$ guarantees that SDP cannot uniquely recover $\bf M^*$. Checking the RIP constant for a linear measurement operator is proven to be NP-hard \citep{tillmann2013computational}. Therefore, it is difficult to solve the problem \eqref{eqn:sdp-inex} analytically. To simplify the problem, we will introduce some notations. We use a matrix representation of the measurement operator $\mathcal A$ as follows:
\[
\mathbf{A} = [\vecc(\mathbf{A}_1), \vecc(\mathbf{A}_2), \dots ,\vecc(\mathbf{A}_d)]^T \in \mathbb{R}^{d \times n^2}.
\]
Then, $\mathbf A\vecc(\mathbf{M})=\mathcal A(\mathbf{M})$ for every matrix $\mathbf{M} \in \RR^{n \times n}$. We define $\mathbf{H} = \mathbf{A}^T\mathbf{A}$, which is the matrix representation of the kernel operator $\mathcal{H} = \mathcal{A}^T\mathcal{A}$ to simplify the last constraint of the problem \eqref{eqn:sdp-inex}.

To derive a RIP bound, we consider the following optimization problem given $\mathbf{M}$ and $\mathbf{M}^*$, where $ \mathbf{M}$ is the global solution of \eqref{eqn:ms-mc-sdp} and $\mathbf{M}^*$ is the ground truth solution:
\begin{equation}\label{eq:ori_lmi}
\begin{aligned}
\min_{\delta,\mathbf H} \quad & \delta \\
\text{s.t.} \quad & \mathbf e^T \mathbf H \mathbf e = 0 \\
& \text{$\mathbf H$ is symmetric and satisfies the $\delta_{2r}$-$\RIP$},
\end{aligned}
\end{equation}
where 
\[
\mathbf e=\vecc(\mathbf{M}^* - \mathbf{M}).
\]
For this fixed $\mathbf{M}$ and $\mathbf{M}^*$, we assume that $\mathbf{M} \neq \mathbf{M}^*$ and that $\rk(\mathbf{M}^*-\mathbf{M}) > 2r$, since if $\rk(\mathbf{M}^*-\mathbf{M}) \leq 2r$, the relation $\mathbf{M}=\mathbf{M}^*$ holds automatically by definition of $\delta_{2r}$-RIP for any $\delta$ since it implies strong convexity. Denote the optimal value to \eqref{eq:ori_lmi} as $\delta(\mathbf e)$, which is a function of $\mathbf e$. It is desirable to find 
\[
	\delta^* \coloneqq \min_{\mathbf e: \tr(\mathbf{M}) \leq \tr(\mathbf{M}^*)} \delta(\mathbf e).
\]
By the logic of in-existence of counterexample,  we know that if a problem $\mathbf{H} = \mathbf{A}^T \mathbf{A}$ has $\delta_{2r}$-RIP with $\delta < \delta^*$, then the solution to \eqref{eqn:ms-mc-sdp} will be $\mathbf{M}^*$, which is the ground truth solution. However, since the last constraint of \eqref{eq:ori_lmi} is non-convex, it is useful to replace it with a surrogate condition that allows solving the problem analytically. The following problem helps to achieve this goal:
\begin{equation}\label{eq:lmi}
\begin{aligned}
\min_{\delta,\mathbf H} \quad & \delta \\
\text{s.t.} \quad & \mathbf e^T \mathbf H \mathbf e \leq 2\|\mathbf{e}_c\|^2+2(l-3) \delta \|\mathbf{e}_c\|^2 \\
& (1-\delta) \mathbf{I}_{n^2} \preceq \mathbf H \preceq (1+\delta)\mathbf{I}_{n^2}.
\end{aligned}
\end{equation}
Here, $l = \lceil n/r \rceil$ and we define $\{\mathbf{e}_i\}_{i=1}^l$ and $\mathbf{e}_c$ in the following fashion. First, consider the eigendecomposition of $\mathbf{M}^*-\mathbf{M}$ and assume that the eigenvalues are ordered in terms of their absolute values, namely, $|\lambda_1| \geq |\lambda_2| \geq \dots \geq |\lambda_n|$. Let $\mathbf{u}_k$'s denote the corresponding orthonormal eigenvectors:
\[
	\mat_S(\mathbf e) = \mathbf{M}^*-\mathbf{M} = \sum_{k=1}^n \lambda_k \mathbf{u}_k \mathbf{u}_k^T.
\]
Then, we define:
\[
	\mathbf{e}_i = \vecc\left(\sum_{k=(i-1)*r+1}^{\min\{i*r,n\}} \lambda_k \mathbf{u}_k \mathbf{u}_k^T\right),
\]
$\mathbf{e}_{2r} = \mathbf{e}_1 + \mathbf{e}_2$, and $\mathbf{e}_{c} = \sum_{i=3}^l \mathbf{e}_i$. The next proposition allows us to replace \eqref{eq:ori_lmi} with \eqref{eq:lmi} because the optimal value of the \eqref{eq:lmi}, $\delta_{lb}(\mathbf{e})$,  gives a lower bound on $\delta(\mathbf{e})$.
\begin{proposition}\label{prop:convex_lb}
    The optimal objective value of the problem \eqref{eq:lmi}, $\delta_{lb}(\mathbf{e})$, is always less than or equal to the optimal objective value of the problem \eqref{eq:ori_lmi},  i.e., $\delta_{\text{lb}}(\mathbf e) \leq \delta(\mathbf e)$. 
\end{proposition}
The proof of this proposition is central to the construction of the sufficiency bound, which is based on using a convex program to serve as an estimate of the non-convex problem. After we extend the RIP$_{2r}$ constraint in \eqref{eq:lmi} to be RIP$_{n}$(thus making it convex), it is necessary to somehow preserve the information that the near isometric property of $\mathbf{H}$ should only apply to low-rank matrices. This is achieved by changing the first constraint so that $\mathbf{e}$ does not need to be completely in the null space of $\mathbf{H}$. \eqref{eq:lmi} approximately requires that $\mathbf{H}$ only maps a certain low-rank sub-manifold to 0. The full proof can be found in the Appendix. As a result of Proposition \ref{prop:convex_lb}, it immediately follows that
\[
	\delta_{\text{lb}} = \min_{\mathbf e: \tr(\mathbf{M}) \leq \tr(\mathbf{M}^*)} \delta_{\text{lb}}(\mathbf e) \leq \delta^*.
\]
In fact, we can obtain a lower bound on the value $\delta_{lb}$ by solving the problem \eqref{eq:lmi} analytically. The following lemma quantifies a lower bound on $\delta_{lb}$.
\begin{lemma}\label{lem:delta_lb}
    It holds that
    \[ 	\delta_{\text{lb}} \geq \frac{2r}{n+(n-2r)(2l-5)}. \]
\end{lemma}
The best-known sufficiency bound presented in \cite{cai2013sharp} is independent of $n$ and $r$. This sufficiency lower bound on the RIP constant presented in Lemma \ref{lem:delta_lb} can be tighter than $1/2$ depending on the size of the problem $n$ and the rank of the ground truth matrix $r$. For instance, the SDP formulation converges to ground truth solution whenever RIP constant $\delta$ is close to $1$ as $r \xrightarrow{} n/2$. On the other hand, whenever $r/n$ is ratio is small, e.g. rank-1 matrix sensing problem with large $n$, $\delta <1/2$ is a stronger guarantee for recovery of the ground truth matrix. Combined with the $1/2$ sufficiency bound that works for both the symmetric and asymmetric cases \cite{cai2013sharp}, we obtain the following result:
\begin{theorem}\label{thm:sdp_rip_new}
    The global solution of the SDP formulation \eqref{eqn:ms-mc-sdp} will be the ground truth matrix $\mathbf{M}^*$ if the sensing matrix $\mathcal{A}$ satisfies the RIP condition with the RIP constant $\delta_{2r}$ satisfying the inequality:
    \[
        \delta_{2r} < \max \left\{1/2, \frac{2r}{n+(n-2r)(2l-5)}\right\},
    \]
    where $l = \lceil n/r \rceil$.
\end{theorem}
Compared with the existing sufficiency RIP bounds, this new result has a striking advantage. The bound $\delta_{2r}<1/2$ has already been proven to be the sharpest for the B-M formulation, which is independent of the search rank. In contrast, Theorem \ref{thm:sdp_rip_new} shows that the RIP bound for SDP exceeds this bound and approaches 1 as the rank $r$ increases. 

In this section, we have shown that as opposed to the popular belief that B-M enjoys very similar RIP guarantees as the SDP approach, there are real benefits to switching to the SDP formulation, making it a more competitive option since specialized SDP solvers are becoming more efficient in recent years. However, we will next provide some problem instances for which the SDP method fails to solve the problem while the B-M method contains no spurious solutions, which balances the desirable properties of the SDP method.

\section{Advantages of the B-M Method}

In this section, we give two classes of rank-1 matrix completion problems for which the B-M factorization does not contain any spurious solution while SDP fails to recover its ground truth matrix. Throughout this section, the rank-1 positive semidefinite ground truth matrix $\mathbf{M}^* = \mathbf{x^*}(\mathbf{x^*})^T$ is assumed not to contain any zero entries, meaning that $x^*_i \not = 0$ for all $i \in [n]$. Before proceeding with the results, we provide two small examples to highlight the underlying ideas behind the main results.
\begin{example}
    Consider a block sparsity graph $\mathcal{G} = (\mathcal{V}, \mathcal{E})$  with $| \mathcal{V}| = 3 $ nodes and the edge set $\mathcal{E} = \{ (1,1), (1,2), (2,3) \}$. Namely, it is a chain graph with 3 nodes and a self-loop at the first node. First, we aim to show that only second-order critical points are the global solutions of the B-M factorization method. The B-M factorization formulation \eqref{eqn:mc-bm} with the squared loss function can be explicitly written as $\min_{ \bf x \in \mathbb{R}^3} f(\mathbf{x})$, where 
    %
    %
    \begin{equation*}
         f(\mathbf{x}) =  \frac{1}{4}\sum_{\substack{(i,j) \in \mathcal{E} \\ i = j}} (x_i^2 - (x^*_i)^2)^2 + \frac{1}{2}\sum_{\substack{(i,j) \in \mathcal{E} \\ i \not = j}} (x_i x_j - x^*_i x^*_j)^2.
    \end{equation*}
    The corresponding gradient and Hessian are:
    %
    %
    
    \begin{align*}
        &\frac{\partial f(\bf x)}{\partial x_i} = \sum_{\substack{i,j \in \mathcal{E} \\ i = j}} (x_i^2 - (x^*_i)^2)x_i + \sum_{\substack{i,j \in \mathcal{E} \\ i \not = j}} (x_i x_j - x^*_i x^*_j)x_j, \\
        &\frac{\partial^2 f(\bf x)}{\partial x_i^2} = \mathbbm{1}[(i,i) \in \mathcal{E}](3 x_i^2 - (x_i^*)^2) + \sum_{ i,j \in \mathcal{E}} x_j^2, \\
        &\frac{\partial^2 f(\bf x)}{\partial x_i \partial x_j} =  \begin{cases}
            2 x_i x_j - x^*_i x^*_j, & \text{if }i \not = j \text{ and } (i,j) \in \mathcal{E} \\
            0, & \text{otherwise}
        \end{cases}.
    \end{align*}

    Each second-order critical point $\hat{\bf x}$ must satisfy the conditions $\nabla f(\hat{ \bf x}) = 0 $ and $\nabla^2 f(\hat{ \bf x}) \succeq 0$. The third entry of the gradient implies either $\hat{x}_2 = 0$ or $\hat{x}_2 \hat{x}_3 = x^*_2 x^*_3$. Whenever $\hat{x}_2 =0$, the Hessian is not positive semidefinite since $[\nabla^2 f(\hat{x})]_{2:3,2:3} \not \succeq 0$. Thus, $\hat{x}_2 \hat{x}_3 = x^*_2 x^*_3$ must hold. Following this, $\partial f(\hat{ \bf x})/ \partial x_2 $ implies either $\hat{x}_1 = 0$ or $\hat{x}_1 \hat{x}_2 = x^*_1 x^*_2$. However, if $\hat{x}_1 = 0$, then $\partial f(\hat{ \bf x})/ \partial x_1 $ gives $ -x^*_1 x^*_2 \hat{x}_2 = 0$, which implies $\hat{x}_2 = 0$. This contradicts the earlier result. Thus, each second-order critical point must have the following properties:
    \[ \hat{x}_1^2  = (x^*_1)^2, \quad \hat{x}_1 \hat{x}_2  = x^*_1 x^*_2, \quad \hat{x}_2 \hat{x}_3  = x^*_2 x^*_3. \]
    The solution to this system of equations proves the exact recovery of the ground truth matrix $\mathbf{M}^*$. Hence, the only second-order critical points are the valid factors of the ground truth solution, i.e $\pm \mathbf{x^*}$. 
    
    The next step is to demonstrate the failure of the SDP formulation \eqref{eqn:ms-mc-sdp} for some instances of the MC problem with this given block sparsity matrix $\mathcal{G}$. The problem \eqref{eqn:ms-mc-sdp} is equivalent to the optimization 
    \begin{align*}
        \min_{\mathbf{M} \in \mathbb{R}^{3 \times 3}} & \mathbf{M}_{2,2} + \mathbf{M}_{3,3} \\
    \text{s.t } &   \begin{bmatrix}
        (x^*_1)^2 & x^*_1 x^*_2 & \mathbf{M}_{1,3} \\ x^*_1 x^*_2 & \mathbf{M}_{2,2} & x^*_2 x^*_3 \\ \mathbf{M}_{3,1} & x^*_2 x^*_3 & \mathbf{M}_{3,3}
    \end{bmatrix} \succeq 0.
    \end{align*}
    Consider a feasible solution $\mathbf{\hat{M}}$ with $\mathbf{\hat{M}}_{2,2} = \mathbf{\hat{M}}_{3,3} = |x^*_2 x^*_3|$ and $\mathbf{\hat{M}}_{1,3} = \mathbf{\hat{M}}_{3,1} = |x^*_1 x^*_2|$. Note that $\mathbf{\hat{M}}$ is feasible whenever $|x^*_3| \geq |x^*_2|$. While the objective value of the ground truth matrix $\mathbf{M}^*$ is $(x^*_2)^2 + (x^*_3)^2$, the objective value of the feasible solution $\mathbf{\hat{M}}$ is $2|x^*_2 x^*_3|$. Under the assumption $|x^*_3| > |x^*_2|$,  the feasible solution $\mathbf{\hat{M}}$ is strictly better than the ground truth solution. Thus, SDP fails to recover the ground truth matrix.
    
\end{example}
This example clearly demonstrates the existence of MC instances for which the B-M method successfully converges to the ground truth solution while the SDP fails to find the solution. One reason is that the number of measurements is $\mathcal{O}(n)$ in this example, which is the minimum threshold for exact completion. However, the statistical guarantees on SDP often need more observations. We can generalize Example 1 to any chain graph with $n$ nodes and a single self-loop at one of the ends.
\begin{theorem} \label{thm: bm-solves-chain}
    Consider the MC problem with a rank-1 positive definite ground truth matrix $\mathbf{M}^* \in \mathbb{R}^{n \times n}$ that can be factorized as $\mathbf{M}^* = \mathbf{x}^* (\mathbf{x}^*)^T$ with $x^*_i \not = 0, \forall i \in [n]$. Let $\mathcal{G} = (\mathcal{V}, \mathcal{E})$ be a block sparsity graph with $|\mathcal{V}| = n$ and $\mathcal{E} = \{ (1,1), (1,2), (2,3), \cdots, (n-1, n) \}$. Then, the B-M method \eqref{eqn:mc-bm} does not contain any spurious solutions.
\end{theorem}
The proof of Theorem \ref{thm: sdp-fails-chain} needs a careful treatment of the second-order optimality conditions and is deferred to the appendix. In addition to the success of the B-M factorization, the next result establishes the failure of the SDP for the instances described in the above theorem. 
\begin{theorem} \label{thm: sdp-fails-chain}
    Consider the ground truth matrix $\mathbf{M}^* \in \mathbb{R}^{n \times n}$ satisfying the conditions in Theorem \ref{thm: sdp-fails-chain}. Suppose that there exist two indices $j,k$ such that $x^*_k > x^*_j$ and $j,k > 2$. Then, the SDP problem \eqref{eqn:ms-mc-sdp} fails to recover the ground truth matrix. 
\end{theorem}
As mentioned before, SDP fails due to a lack of observations on the diagonal entries of the ground truth matrix. Note that the RIP condition is not satisfied since these are MC problems. As a result, whenever we do not have sufficient guarantees on linear measurement operator, none of the methods are superior to the other one in terms of exact recovery. The next example identifies another class of problem instances that corroborates these findings. 
\begin{example}
    Consider a block sparsity graph $\mathcal{G} = (\mathcal{V}, \mathcal{E})$  with $| \mathcal{V}| = 3 $ nodes and the edge set $\mathcal{E} = \{ (1,2), (2,3), (3,1) \}$. Namely, it is a simple cycle with 3 nodes. The B-M factorization formulation \eqref{eqn:mc-bm} with the squared loss function can be written the same as in Example 1. Firstly, we can show that each second-order critical point $\mathbf{\hat{x}}$ only has nonzero entries, i.e., $\hat{x}_i \not = 0$ for all $i \in [n]$. Without loss of generality, suppose by contradiction that $\hat{x}_1 = 0$. In order for the stationarity condition to hold, either $\hat{x}_2 \hat{x}_3 = x^*_2 x^*_3$ or $\hat{x}_2 = \hat{x}_3 = 0$ should be satisfied. The latter implies $\mathbf{\hat{x}} = 0$ and $\nabla^2 f(\mathbf{\hat{x}}) \not \succeq 0$ in that case. Thus, $\hat{x}_2 \hat{x}_3 = x^*_2 x^*_3$ must hold. Following this, $\partial f(\mathbf{\hat{x}})/ \partial x_1$ yields $-x^*_1 x^*_2 \hat{x}_2 - x^*_1 x^*_3 \hat{x}_3 =0$. Combining these two equations yields $(\hat{x}_3)^2 = -(x^*_2)^2$, which does not have any real solution. As a result, each second-order critical point $\mathbf{\hat{x}}$ must have only nonzero entries.
    %
    %
    %
    %
    %
    
    By the condition $\nabla f(\mathbf{\hat{x}}) = 0$, whenever $\hat{x}_i \hat{x}_j = x^*_i x^*_j$ for some $(i,j) \in \mathcal{E}$, then $\hat{x}_i \hat{x}_j = x^*_i x^*_j$ holds for every $(i,j) \in \mathcal{E}$. This system of equations yields the ground truth solution. Accordingly, a spurious solution $\mathbf{\hat{x}}$ must have the following characteristics: $\hat{x}_i \hat{x}_j \not = x^*_i x^*_j, \forall (i,j) \in \mathcal{E}$ and $\hat{x}_i \not = 0, \forall i \in \{1,2,3 \} $. Define $a_{i,j} = x_i x_j - x^*_i x^*_j  $. Then, the stationarity condition becomes
    \begin{equation*}
        \nabla f(\mathbf{\hat{x}}) = \begin{bmatrix}
            \hat{a}_{1,2}\hat{x}_2 +  \hat{a}_{1,3}\hat{x}_3 \\
            \hat{a}_{1,2}\hat{x}_1 + \hat{a}_{2,3}\hat{x}_3 \\
            \hat{a}_{1,3}\hat{x}_1 + \hat{a}_{2,3}\hat{x}_2
        \end{bmatrix} = 0.
    \end{equation*}
    Multiplying the first entry of the gradient by $\hat{x}_1 / \hat{x}_2$ and substituting with the second entry gives $\hat{a}_{1,2} \hat{x}_1 = - \hat{a}_{2,3} \hat{x}_3$. Successively, substituting this to the third entry of the gradient results in 
    \[ \hat{x}_3 (\hat{a}_{1,3} \hat{x}_1 - \hat{a}_{2,3} \hat{x}_2) = 0. \]
    Because we search for a solution with $\hat{x}_i \not = 0$, we must have $\hat{a}_{1,3} \hat{x}_1 - \hat{a}_{2,3} \hat{x}_2 = 0$. This condition combined with the last entry of the gradient results in the condition $\hat{a}_{1,3} \hat{x}_1 = 0$, which is a contradiction. As a result, all the second-order critical points are global solutions that yield the ground truth matrix completion.
    
     Our next goal is to show that the SDP formulation \eqref{eqn:ms-mc-sdp} fails for this class of instances of MC instances. Note that the SDP formulation of the matrix completion problem considered in this example is equivalent to the formulation: 
    \begin{align*}
        \min_{\mathbf{M} \in \mathbb{R}^{3 \times 3}} & \mathbf{M}_{1,1} + \mathbf{M}_{2,2} + \mathbf{M}_{3,3} \\
    \text{s.t } &   \begin{bmatrix}
        \mathbf{M}_{1,1} & x^*_1 x^*_2 & x^*_1 x^*_3 \\ x^*_1 x^*_2 & \mathbf{M}_{2,2} & x^*_2 x^*_3 \\ x^*_1 x^*_3 & x^*_2 x^*_3 & \mathbf{M}_{3,3}
    \end{bmatrix} \succeq 0.
    \end{align*}
    Without loss of generality, assume that $x^*_1 \leq x^*_2 \leq x^*_3$ by the symmetry of the problem. Consider a feasible rank-2 solution $\mathbf{\hat{M}}$ given as $\mathbf{\hat{M}}_{1,1} = x^*_1 (x^*_3 - x^*_2), \mathbf{\hat{M}}_{2,2} = x^*_2 (x^*_3 - x^*_1) $ and $\mathbf{M}_{3,3} = x^*_3 (x^*_1 + x^*_2)$. Note that $\mathbf{\hat{M}}$ is feasible whenever $x^*_3 \geq x^*_1 + x^*_2$. The objective value of the feasible solution $\mathbf{\hat{M}}$ is $ -2 x^*_1 x^*_2 + 2 x^*_2 x^*_3 + 2 x^*_1 x^*_3$,  whereas the objective of the ground truth solution $\mathbf{M}^*$ is $(x^*_1)^2 + (x^*_2)^2 + (x^*_3)^2$. The feasible solution $\mathbf{\hat{M}}$ is strictly better than the ground truth solution if $x^*_3 > x^*_1 + x^*_2$. Hence, the SDP cannot recover the ground truth solution for all the instances with a simple cycle block sparsity graph.
\end{example}
Similar to Example 1, SDP fails in this example due to a lack of diagonal observations. Next, we can generalize this instance to any simple cycle block sparsity graph with an odd number of vertices. 
\begin{theorem} \label{thm: bm-solves-cycle}
    Consider the matrix completion problem with a rank-1 positive definite ground truth matrix $\mathbf{M}^* \in \mathbb{R}^{n \times n}$ that can be factorized as $\mathbf{M}^* = \mathbf{x}^* (\mathbf{x}^*)^T$ with $x^*_i \not = 0, \forall i \in [n]$. Let $\mathcal{G} = (\mathcal{V}, \mathcal{E})$ be a block sparsity graph with $|\mathcal{V}| = |\mathcal{E}| = n = 2k+1$ and $\mathcal{E} = \{ (0,1), (1,2), \dots (2k-1, 2k), (2k, 0) \}$. Then, the B-M factorization problem \eqref{eqn:mc-bm} does not contain any spurious solutions. 
\end{theorem} 
\begin{theorem} \label{thm: sdp-fails-cycle}
    Consider the ground truth matrix $\mathbf{M}^* \in \mathbb{R}^{n \times n}$ satisfying the conditions in Theorem \ref{thm: bm-solves-chain}. Suppose that the condition $\sum_{t=1}^{k} (x^*_{2t-1})^2 > \sum_{t=0}^{k} (x^*_{2t})^2$ holds for some chosen node $0$. Then, the SDP problem \eqref{eqn:ms-mc-sdp} fails to recover the ground truth matrix.
\end{theorem}
Note that we can choose any node as node $0$ due to the symmetry of the problem. Therefore, the condition stated in Theorem \ref{thm: sdp-fails-cycle} is not restrictive because this condition suffices to hold for a chosen node $0$ among $2k+1$ ones.
As a result of the above theorems, the B-M factorization can outperform the convex relaxation approach. One important extension of the work presented in this section would be finding subgraphs for which the B-M method is successful while the SDP is unsuccessful, and then attaching these subgraphs to generate larger block sparsity graphs. It is known that SDP will fail for these instances and it is intriguing to investigate the behavior of the B-M method for those instances.

\section{Conclusions}

In this paper, we conducted a comparison between two main approaches to the matrix completion and matrix sensing problems: a convex relaxation that gives an SDP formulation and the B-M factorization method. It is well-known that both of these methods enjoy mathematical guarantees for the recovery of the ground truth matrix whenever the RIP assumption is satisfied with a sufficiently small $\delta$. We offered the first result in the literature that compares these two methods whenever the RIP condition is not satisfied or only satisfied with a large constant. We discovered classes of problems for which B-M factorization fails while the SDP recovers the ground truth matrix. The fact that specialized SDP algorithms are improved in recent years and can compete with simple first-order descent algorithms inspired us to investigate sharper bounds on sufficient conditions for the SDP formulation. We provided RIP bounds for the SDP formulation that depend on the rank of the solution and are automatically satisfied for high-rank problems, unlike the B-M method. On the other hand, when the number of measurements from the ground truth matrix is not high, we showed that SDP fails drastically while the B-M method does not contain any spurious solutions on its optimization landscape. As a result, we conclude that none of the methods outperforms the other one whenever the sufficiency guarantees are not met. The parameters of the problem, such as dimension, rank, and linear measurement operator, determine which solution method performs better. Consequently, it is prudent to apply both solution methods in case the RIP and incoherence are not satisfied. 

\bibliography{main}

\begin{thebibliography}{39}
\providecommand{\natexlab}[1]{#1}

\bibitem[{Bhojanapalli, Neyshabur, and Srebro(2016)}]{bhojanapalli2016global}
Bhojanapalli, S.; Neyshabur, B.; and Srebro, N. 2016.
\newblock Glob\-al Optimality of Local Search for Low Rank Matrix Recovery.
\newblock In \emph{Advances in Neural Information Processing Systems},
  volume~29.

\bibitem[{Bi and Lavaei(2021)}]{bi2020global}
Bi, Y.; and Lavaei, J. 2021.
\newblock On the absence of spurious local minima in nonlinear low-rank matrix
  recovery problems.
\newblock In \emph{International Conference on Artificial Intelligence and
  Statistics}, 379--387. PMLR.

\bibitem[{Boumal(2016)}]{boumal2016nonconvex}
Boumal, N. 2016.
\newblock Nonconvex phase synchronization.
\newblock \emph{SIAM Journal on Optimization}, 26(4): 2355--2377.

\bibitem[{Burer and Monteiro(2003)}]{burer2003nonlinear}
Burer, S.; and Monteiro, R.~D. 2003.
\newblock A nonlinear programming algorithm for solving semidefinite programs
  via low-rank factorization.
\newblock \emph{Mathematical Programming}, 95(2): 329--357.

\bibitem[{Cai and Zhang(2013)}]{cai2013sharp}
Cai, T.~T.; and Zhang, A. 2013.
\newblock Sharp RIP bound for sparse signal and low-rank matrix recovery.
\newblock \emph{Applied and Computational Harmonic Analysis}, 35(1): 74--93.

\bibitem[{Candes et~al.(2015)Candes, Eldar, Strohmer, and
  Voroninski}]{candes2015phase}
Candes, E.~J.; Eldar, Y.~C.; Strohmer, T.; and Voroninski, V. 2015.
\newblock Phase retrieval via matrix completion.
\newblock \emph{SIAM review}, 57(2): 225--251.

\bibitem[{Candes and Plan(2010)}]{candes2010tight}
Candes, E.~J.; and Plan, Y. 2010.
\newblock Tight oracle bounds for low-rank matrix recovery from a minimal
  number of random measurements.
\newblock \emph{arXiv preprint arXiv:1001.0339}.

\bibitem[{Cand{\`e}s and Recht(2009)}]{candes2009exact}
Cand{\`e}s, E.~J.; and Recht, B. 2009.
\newblock Exact matrix completion via convex optimization.
\newblock \emph{Foundations of Computational mathematics}, 9(6): 717--772.

\bibitem[{Cand{\`e}s and Tao(2010)}]{candes2010power}
Cand{\`e}s, E.~J.; and Tao, T. 2010.
\newblock The power of convex relaxation: Near-optimal matrix completion.
\newblock \emph{IEEE Transactions on Information Theory}, 56(5): 2053--2080.

\bibitem[{Davenport et~al.(2014)Davenport, Plan, Van Den~Berg, and
  Wootters}]{davenport20141}
Davenport, M.~A.; Plan, Y.; Van Den~Berg, E.; and Wootters, M. 2014.
\newblock 1-bit matrix completion.
\newblock \emph{Information and Inference: A Journal of the IMA}, 3(3):
  189--223.

\bibitem[{Fattahi and Sojoudi(2020)}]{fattahi2020exact}
Fattahi, S.; and Sojoudi, S. 2020.
\newblock Exact guarantees on the absence of spurious local minima for
  non-negative rank-1 robust principal component analysis.
\newblock \emph{Journal of machine learning research}.

\bibitem[{Ge, Jin, and Zheng(2017)}]{ge2017no}
Ge, R.; Jin, C.; and Zheng, Y. 2017.
\newblock No spurious local minima in nonconvex low rank problems: A unified
  geometric analysis.
\newblock In \emph{International Conference on Machine Learning}, 1233--1242.
  PMLR.

\bibitem[{Ha, Liu, and Barber(2020)}]{HLB2020}
Ha, W.; Liu, H.; and Barber, R.~F. 2020.
\newblock An Equivalence Between Critical Points for Rank Constraints Versus
  Low-Rank Factorizations.
\newblock \emph{SIAM Journal on Optimization}, 30(4): 2927--2955.

\bibitem[{Jin et~al.(2020)Jin, Lavaei, Sojoudi, and Baldick}]{jin2020boundary}
Jin, M.; Lavaei, J.; Sojoudi, S.; and Baldick, R. 2020.
\newblock Boundary defense against cyber threat for power system state
  estimation.
\newblock \emph{IEEE Transactions on Information Forensics and Security}, 16:
  1752--1767.

\bibitem[{Jin et~al.(2019)Jin, Molybog, Mohammadi-Ghazi, and
  Lavaei}]{jin2019towards}
Jin, M.; Molybog, I.; Mohammadi-Ghazi, R.; and Lavaei, J. 2019.
\newblock Towards robust and scalable power system state estimation.
\newblock In \emph{2019 IEEE 58th Conference on Decision and Control (CDC)},
  3245--3252. IEEE.

\bibitem[{Koren, Bell, and Volinsky(2009)}]{koren2009matrix}
Koren, Y.; Bell, R.; and Volinsky, C. 2009.
\newblock Matrix factorization techniques for recommender systems.
\newblock \emph{Computer}, 42(8): 30--37.

\bibitem[{Ma et~al.(2022)Ma, Bi, Lavaei, and Sojoudi}]{ma2021sharp}
Ma, Z.; Bi, Y.; Lavaei, J.; and Sojoudi, S. 2022.
\newblock Sharp Restricted Isometry Property Bounds for Low-rank Matrix
  Recovery Problems with Corrupted Measurements.
\newblock \emph{AAAI-22}.

\bibitem[{Ma and Sojoudi(2022)}]{ma2022noisy}
Ma, Z.; and Sojoudi, S. 2022.
\newblock Noisy Low-rank Matrix Optimization: Geometry of Local Minima and
  Convergence Rate.
\newblock \emph{arXiv preprint arXiv:2203.03899}.

\bibitem[{Mohan and Fazel(2010)}]{mohan2010new}
Mohan, K.; and Fazel, M. 2010.
\newblock New restricted isometry results for noisy low-rank recovery.
\newblock In \emph{2010 IEEE International Symposium on Information Theory},
  1573--1577. IEEE.

\bibitem[{Park et~al.(2017)Park, Kyrillidis, Carmanis, and Sanghavi}]{PKCS2017}
Park, D.; Kyrillidis, A.; Carmanis, C.; and Sanghavi, S. 2017.
\newblock Non-square Matrix Sensing Without Spurious Local Minima via the
  {Burer}--{Monteiro} Approach.
\newblock In \emph{Proceedings of the 20th International Conference on
  Artificial Intelligence and Statistics}, volume~54 of \emph{Proceedings of
  Machine Learning Research}, 65--74.

\bibitem[{Recht, Fazel, and Parrilo(2010)}]{recht2010guaranteed}
Recht, B.; Fazel, M.; and Parrilo, P.~A. 2010.
\newblock Guaranteed minimum-rank solutions of linear matrix equations via
  nuclear norm minimization.
\newblock \emph{SIAM review}, 52(3): 471--501.

\bibitem[{Recht, Xu, and Hassibi(2008)}]{recht2008necessary}
Recht, B.; Xu, W.; and Hassibi, B. 2008.
\newblock Necessary and sufficient conditions for success of the nuclear norm
  heuristic for rank minimization.
\newblock In \emph{2008 47th IEEE Conference on Decision and Control},
  3065--3070. IEEE.

\bibitem[{Shechtman et~al.(2015)Shechtman, Eldar, Cohen, Chapman, Miao, and
  Segev}]{shechtman2015phase}
Shechtman, Y.; Eldar, Y.~C.; Cohen, O.; Chapman, H.~N.; Miao, J.; and Segev, M.
  2015.
\newblock Phase retrieval with application to optical imaging: a contemporary
  overview.
\newblock \emph{IEEE signal processing magazine}, 32(3): 87--109.

\bibitem[{Singer(2011)}]{singer2011angular}
Singer, A. 2011.
\newblock Angular synchronization by eigenvectors and semidefinite programming.
\newblock \emph{Applied and Computational Harmonic Analysis}, 30(1): 20--36.

\bibitem[{Tillmann and Pfetsch(2013)}]{tillmann2013computational}
Tillmann, A.~M.; and Pfetsch, M.~E. 2013.
\newblock The computational complexity of the restricted isometry property, the
  nullspace property, and related concepts in compressed sensing.
\newblock \emph{IEEE Transactions on Information Theory}, 60(2): 1248--1259.

\bibitem[{Yal{\c{c}}{\i}n et~al.(2022)Yal{\c{c}}{\i}n, Zhang, Lavaei, and
  Sojoudi}]{yalccin2022factorization}
Yal{\c{c}}{\i}n, B.; Zhang, H.; Lavaei, J.; and Sojoudi, S. 2022.
\newblock Factorization approach for low-complexity matrix completion problems:
  Exponential number of spurious solutions and failure of gradient methods.
\newblock In \emph{International Conference on Artificial Intelligence and
  Statistics}, 319--341. PMLR.

\bibitem[{Yurtsever et~al.(2021)Yurtsever, Tropp, Fercoq, Udell, and
  Cevher}]{yurtsever21}
Yurtsever, A.; Tropp, J.~A.; Fercoq, O.; Udell, M.; and Cevher, V. 2021.
\newblock Scalable Semidefinite Programming.
\newblock \emph{SIAM Journal on Mathematics of Data Science}, 3(1): 171--200.

\bibitem[{Yurtsever et~al.(2017)Yurtsever, Udell, Tropp, and
  Cevher}]{udell2017}
Yurtsever, A.; Udell, M.; Tropp, J.; and Cevher, V. 2017.
\newblock {Sketchy Decisions: Convex Low-Rank Matrix Optimization with Optimal
  Storage}.
\newblock In Singh, A.; and Zhu, J., eds., \emph{Proceedings of the 20th
  International Conference on Artificial Intelligence and Statistics},
  volume~54 of \emph{Proceedings of Machine Learning Research}, 1188--1196.
  PMLR.

\bibitem[{Zhang and Zhang(2020)}]{zhang2020many}
Zhang, G.; and Zhang, R.~Y. 2020.
\newblock How Many Samples Is a Good Initial Point Worth in Low-Rank Matrix
  Recovery?
\newblock In \emph{Advances in Neural Information Processing Systems},
  volume~33, 12583--12592.

\bibitem[{Zhang, Bi, and Lavaei(2021{\natexlab{a}})}]{zhang2021general}
Zhang, H.; Bi, Y.; and Lavaei, J. 2021{\natexlab{a}}.
\newblock General low-rank matrix optimization: Geometric analysis and sharper
  bounds.
\newblock \emph{arXiv preprint arXiv:2104.10356}.

\bibitem[{Zhang, Bi, and Lavaei(2021{\natexlab{b}})}]{ZBL2021}
Zhang, H.; Bi, Y.; and Lavaei, J. 2021{\natexlab{b}}.
\newblock General Low-Rank Matrix Optimization: Geometric Analysis and Sharper
  Bounds.
\newblock In \emph{Advances in Neural Information Processing Systems}.

\bibitem[{Zhang et~al.(2022)Zhang, Yalcin, Lavaei, and
  Sojoudi}]{zhang2022unified}
Zhang, H.; Yalcin, B.; Lavaei, J.; and Sojoudi, S. 2022.
\newblock A Unified Complexity Metric for Nonconvex Matrix Completion and
  Matrix Sensing in the Rank-one Case.
\newblock \emph{arXiv preprint arXiv:2204.02364}.

\bibitem[{Zhang(2021)}]{Zhang2021-p}
Zhang, R.~Y. 2021.
\newblock Sharp Global Guarantees for Nonconvex Low-Rank Matrix Recovery in the
  Overparameterized Regime.
\newblock ArXiv:2104.10790.

\bibitem[{Zhang and Lavaei(2021)}]{zhang2021sparse}
Zhang, R.~Y.; and Lavaei, J. 2021.
\newblock Sparse semidefinite programs with guaranteed near-linear time
  complexity via dualized clique tree conversion.
\newblock \emph{Mathematical programming}, 188(1): 351--393.

\bibitem[{Zhang, Sojoudi, and Lavaei(2019)}]{zhang2019sharp}
Zhang, R.~Y.; Sojoudi, S.; and Lavaei, J. 2019.
\newblock Sharp Restricted Isometry Bounds for the Inexistence of Spurious
  Local Minima in Nonconvex Matrix Recovery.
\newblock \emph{Journal of Machine Learning Research}, 20(114): 1--34.

\bibitem[{Zhang et~al.(2018)Zhang, Wang, Yu, and Gu}]{ZWYG2018}
Zhang, X.; Wang, L.; Yu, Y.; and Gu, Q. 2018.
\newblock A Primal-Dual Analysis of Global Optimality in Nonconvex Low-Rank
  Matrix Recovery.
\newblock In \emph{Proceedings of the 35th International Conference on Machine
  Learning}, volume~80 of \emph{Proceedings of Machine Learning Research},
  5862--5871.

\bibitem[{Zhang, Madani, and Lavaei(2017)}]{zhang2017conic}
Zhang, Y.; Madani, R.; and Lavaei, J. 2017.
\newblock Conic relaxations for power system state estimation with line
  measurements.
\newblock \emph{IEEE Transactions on Control of Network Systems}, 5(3):
  1193--1205.

\bibitem[{Zhu et~al.(2018)Zhu, Li, Tang, and Wakin}]{ZLTW2018}
Zhu, Z.; Li, Q.; Tang, G.; and Wakin, M.~B. 2018.
\newblock Global Optimality in Low-Rank Matrix Optimization.
\newblock \emph{IEEE Transactions on Signal Processing}, 66(13): 3614--3628.

\bibitem[{Zhu et~al.(2021)Zhu, Li, Tang, and Wakin}]{zhu2021global}
Zhu, Z.; Li, Q.; Tang, G.; and Wakin, M.~B. 2021.
\newblock The global optimization geometry of low-rank matrix optimization.
\newblock \emph{IEEE Transactions on Information Theory}, 67(2): 1308--1331.

\end{thebibliography}

\appendix
\onecolumn

\section{Appendices}

\subsection{Proof of Theorem \ref{thm: sdp-solves-mc}}

\begin{proof}
    Let $\mathbf{M}^* \in \mathbb{R}^{n \times n}$ be the rank-r unperturbed ground truth $m \times m$ block matrix with each block having dimension $r \times r$. Hence, the ground truth matrix can be factorized as $\mathbf{M}^* = \mathbf{X}^* (\mathbf{X}^*)^T, \mathbf{X}^* \in \mathbb{R}^{n \times r}$. Each square factor $\mathbf{X}^*_i \in \mathbb{R}^{r \times r}$ is rank-r if $i \in \mathcal{S}_1$ and is $\mathbf{0}_{r \times r}$ otherwise. We perturb the ground truth matrix by $\epsilon \in \mathbb{R}^{n \times r}$, where $\mathbf{M}^*(\epsilon) = (\mathbf{X}^* + \epsilon)(\mathbf{X}^* + \epsilon)^T$ such that $\mathbf{X}^*(\epsilon)_i = \mathbf{X}^*_i + \epsilon_i$ if $i \in \mathcal{S}_1$ and $\mathbf{X}^*(\epsilon)_i = \epsilon_i$ otherwise. Here, we assume that $\rk(\mathbf{X}^*_i + \epsilon_i) = r, \forall i \in [m]$. 
    
    If all diagonal blocks are observed, then it will reduce to a feasibility problem and we can skip the following procedure. Since $\mathcal{S}_1$ is maximal independent set, there exists two indices $i \in \mathcal{S}_1$ and $j \not \in \mathcal{S}_1$ with $(j,j) \not \in \mathcal{E}$ such that $(i,j) \in \mathcal{E}_1$. Consider the $2 \times 2$ block sub-matrix with $i$-th and $j$-th block columns and rows:
     \begin{align*}
        \begin{bmatrix} \mathbf{M}_{i,i} & \mathbf{M}_{i,j} \\ \mathbf{M}_{j,i} & \mathbf{M}_{j,j} \end{bmatrix} = \begin{bmatrix}
            (\mathbf{X}^*_i + \epsilon_i)(\mathbf{X}^*_i + \epsilon_i)^T & (\mathbf{X}^*_i + \epsilon_i)\epsilon_j^T \\ \epsilon_j(\mathbf{X}^*_i + \epsilon_i)^T & \mathbf{M}_{j,j}
        \end{bmatrix} \succeq 0.
    \end{align*}
    The equality holds because the blocks $(i,i)$ and $(i,j)$ are in $\mathcal{E}_1$. By the Schur complement argument and since $\mathbf{X}^*_i + \epsilon_i$ is full-rank, the above constraint is equivalent to
    \begin{align*}
        \mathbf{M}_{j,j} \succeq \epsilon_j \epsilon_j^T.
    \end{align*}
    The unique trace minimizer for the diagonal blocks is $\mathbf{M}_{j,j} = \epsilon_j \epsilon_j^T$. By the same argument, this must hold for every $j$ that is not in the independent set without a self-loop. Thus, the objective value cannot be less than $\sum_{i=1}^n \mathrm{tr}(\mathbf{X}^*(\epsilon)_i (\mathbf{X}^*(\epsilon)_i)^T)$. Therefore, the minimum value is achieved whenever $\mathbf{M}^*_{i,i} =  \mathbf{X}^*(\epsilon)_i (\mathbf{X}^*(\epsilon))^T$. This makes $\mathbf{M}(\epsilon)^*$ an optimal solution.
    
    We now prove the uniqueness of the solution. Since the graph is connected, there exists a node $k$ adjacent to the node $j$ such that the edges $(i,j)$ and $(j,k)$ exist in the graph. Consider the $3 \times 3$ block submatrix 
    \begin{align*}
        \begin{bmatrix} \mathbf{M}_{i,i} & \mathbf{M}_{i,j} & \mathbf{M}_{i,k}  \\ \mathbf{M}_{j,i} & \mathbf{M}_{j,j} & \mathbf{M}_{j,k} \\ \mathbf{M}_{k,i} & \mathbf{M}_{k,j} & \mathbf{M}_{k,k}  \end{bmatrix} = \begin{bmatrix} \mathbf{X}^*(\epsilon)_i \mathbf{X}^*(\epsilon)_i^T & \mathbf{X}^*(\epsilon)_i \mathbf{X}^*(\epsilon)_j^T & \mathbf{M}_{i,k}  \\ \mathbf{X}^*(\epsilon)_j \mathbf{X}^*(\epsilon)_i^T & \mathbf{X}^*(\epsilon)_j \mathbf{X}^*(\epsilon)_j^T & \mathbf{X}^*(\epsilon)_j \mathbf{X}^*(\epsilon)_k^T \\ \mathbf{M}_{k,i} & \mathbf{X}^*(\epsilon)_k \mathbf{X}^*(\epsilon)_j^T & \mathbf{X}^*(\epsilon)_k\mathbf{X}^*(\epsilon)_k^T   \end{bmatrix} \succeq 0.
    \end{align*}
    This is equivalent to following constraints by Schur complement
    \begin{align*}
        \begin{bmatrix} \mathbf{X}^*(\epsilon)_i \mathbf{X}^*(\epsilon)_i^T & \mathbf{X}^*(\epsilon)_i \mathbf{X}^*(\epsilon)_j^T  \\ \mathbf{X}^*(\epsilon)_j \mathbf{X}^*(\epsilon)_i^T & \mathbf{X}^*(\epsilon)_j \mathbf{X}^*(\epsilon)_j^T  \end{bmatrix} - \begin{bmatrix}
            \mathbf{M}_{i,k}  \mathbf{X}^*(\epsilon)_k^{-T} \mathbf{X}^*(\epsilon)_k^{-1} \mathbf{M}_{k,i} & \mathbf{M}_{i,k}\mathbf{X}^*(\epsilon)_k^{-T}\mathbf{X}^*(\epsilon)_j^{T} \\
            \mathbf{X}^*(\epsilon)_j \mathbf{X}^*(\epsilon)_k^{-1} \mathbf{M}_{k,i} & \mathbf{X}^*(\epsilon)_j\mathbf{X}^*(\epsilon)_j^{T}
        \end{bmatrix} & \succeq 0, \\
        \begin{bmatrix}
            \mathbf{X}^*(\epsilon)_i\mathbf{X}^*(\epsilon)_i^T - \mathbf{M}_{i,k}  \mathbf{X}^*(\epsilon)_k^{-T}\mathbf{X}^*(\epsilon)_k^{-1} \mathbf{M}_{k,i} & (\mathbf{X}^*(\epsilon)_i - \mathbf{M}_{i,k}\mathbf{X}^*(\epsilon)_k^{-T} )\mathbf{X}^*(\epsilon)_j^T \\
            \mathbf{X}^*(\epsilon)_j(\mathbf{X}^*(\epsilon)_i - \mathbf{M}_{i,k}\mathbf{X}^*(\epsilon)_k^{-T} )^T & 0
        \end{bmatrix} & \succeq 0.
    \end{align*}
    Another Schur complement argument gives 
    \begin{align*}
        (\mathbf{X}^*(\epsilon)_i - \mathbf{M}_{i,k}\mathbf{X}^*(\epsilon)_k^{-T} )\mathbf{X}^*(\epsilon)_j^T = 0.
    \end{align*}
    Since $\mathbf{X}^*(\epsilon)_j^T$ is full rank, we have $\mathbf{X}^*(\epsilon)_i - \mathbf{M}_{i,k}\mathbf{X}^*(\epsilon)_k^{-T} = 0 $. Thus, we obtain $\mathbf{M}_{i,k} = \mathbf{X}^*(\epsilon)_i \mathbf{X}^*(\epsilon)_k^T$. Note that filling out the unobserved non-diagonal blocks is equivalent to adding the edge $(i,k)$ to the graph $\mathcal{G}_1$. Hence, we can always find such triple $(i,j,k)$ defined as above until filling out all missing entries. As a result, we obtain the unique solution $\mathbf{M}^*(\epsilon)$ by continuing iteratively.
\end{proof}

\subsection{Proof of Proposition \ref{prop:convex_lb}}

\begin{proof}
     To prove Proposition \ref{prop:convex_lb}, we study intermediary problem.
    \begin{equation}\label{eq:int_lmi}
    \begin{aligned}
    \min_{\delta,\hat{\mathbf H}} \quad & \delta \\
    \text{s.t.} \quad & \hat{\mathbf e}^T \hat{\mathbf H} \hat{\mathbf e} \leq (1+\delta)\|\mathbf{e}_c\|^2+2(l-3) \delta \|\mathbf{e}_c\|^2 \\
    & (1-\delta)\mathbf{I}_{4r^2} \preceq \hat{\mathbf H} \preceq (1+\delta)\mathbf{I}_{4r^2},
    \end{aligned}
    \end{equation}
    where
    \[
    	\hat{\mathbf e} = \mathbf{P}^T \mathbf{e}, \qquad \mathbf{P} \in \RR^{n^2 \times 4r^2}= P \otimes P,
    \]
    and $P \in \RR^{n \times 2r}$ is defined to be
    \[
    	P = \begin{bmatrix}
    		\mathbf{u}_1 & \mathbf{u}_2 & \dots & \mathbf{u}_{2r},
    	\end{bmatrix}
    \]
    where $\mathbf{u}_i$'s are orthonormal eigenvectors of $\mathbf{M}^*-\mathbf{M}$ so that $P^T P = \mathbf{I}$. Denote the optimal solution to \eqref{eq:int_lmi} as $\delta_P(\mathbf e)$. Then, the following two lemmas will suffice to prove Proposition \ref{prop:convex_lb}.
    \begin{lemma}\label{lem:delta_ineq_1}
    	Given a fixed vector $\mathbf e \in \RR^{n^2}$, we have
    	\begin{equation}
    		\delta_P(\mathbf e) \leq \delta(\mathbf e).
    	\end{equation}
    \end{lemma}
    \begin{lemma}\label{lem:delta_ineq_2}
    	Given a fixed vector $\mathbf e \in \RR^{n^2}$, we have
    	\begin{equation}
    		\delta_{\text{lb}}(\mathbf e) \leq \delta_P(\mathbf e).
    	\end{equation}
    \end{lemma}
    \begin{proof}[Proof of Lemma \ref{lem:delta_ineq_1}]
    	It suffices to show that for any feasible pair $(\delta, \bar{\mathbf H})$ of \eqref{eq:ori_lmi}, we can construct a feasible solution $(\delta,\hat{\mathbf H})$ to \eqref{eq:int_lmi} characterized as below
    	\[
    		\delta = \delta, \qquad \hat{\mathbf H} = \mathbf{P}^T \bar{\mathbf H} \mathbf{P},
    	\]
    	 which directly proves the lemma. We can verify the feasibility of $(\delta,\hat{\mathbf H})$ as follows. The feasibility of the first constraint is certified by the following argument:
    	 \begin{align*}
    	 	\hat{\mathbf e}^T \hat{\mathbf H} \hat{\mathbf e} = \mathbf{e}^T \mathbf{P} \mathbf{P}^T \bar{\mathbf H} \mathbf{P} \mathbf{P}^T \mathbf e,
    	 \end{align*}
    	 By the definition of $\mathbf{P}$, one can write
    	\[
    		\mathbf{P} \mathbf{P}^T \mathbf e = (PP^T \otimes PP^T) \mathbf e = \vecc(PP^T (\mathbf{M}^*-\mathbf{M}) PP^T) = \mathbf{e}_1 + \mathbf{e}_2 = \mathbf{e}_{2r},
    	\]
    	Since $\mathbf e^T \bar{\mathbf H} \mathbf e = 0$ and $\bar{\mathbf H}$ is symmetric,  $\bar{\mathbf H}$ admits a factorization $\bar{\mathbf H} = \bar {\mathbf A}^T \bar{ \mathbf A}$, making $ \bar{ \mathbf A} \mathbf e = 0$. Also, we know that $\mathbf e = \mathbf e_{2r} + \mathbf e_c$, meaning that
    	\[
    		\bar{ \mathbf A} \mathbf e_{2r} = -\bar{ \mathbf A} \mathbf e_{c}.
    	\]
    	Therefore,
    	\begin{align*}
    		\hat{\mathbf e}^T \hat{\mathbf H} \hat{\mathbf e} = \mathbf e_{2r}^T \bar{\mathbf H} \mathbf e_{2r} = \mathbf e_{c}^T \bar{\mathbf H} \mathbf e_{c} =(\sum_{i=3}^l \mathbf{e}_i)^T \bar{\mathbf H} (\sum_{i=3}^l \mathbf{e}_i).
    	\end{align*}
    	Since $\bar{\mathbf H}$ satisfies $\delta_{2r}$-RIP, for every $(i,j)$ such that $i \neq j$, we have:
    	\begin{equation}\label{eq:ec_ineq_1}
    		(\mathbf{e}_i+\mathbf{e}_j)^T \bar{\mathbf H} (\mathbf{e}_i+\mathbf{e}_j) \leq (1+\delta) \|\mathbf{e}_i+\mathbf{e}_j\|^2 = (1+\delta)(\|\mathbf{e}_i\|^2 + \|\mathbf{e}_j\|^2),
    	\end{equation}
    	where the last equality follows from the facts that  $\mathbf{e}_i^T \mathbf{e}_j = 0$ and
    	\begin{equation}\label{eq:ec_ineq_2}
    		(\mathbf{e}_i+\mathbf{e}_j)^T \bar{\mathbf H} (\mathbf{e}_i+\mathbf{e}_j) = \mathbf{e}_i^T \bar{\mathbf H} \mathbf{e}_i + 2\mathbf{e}_i^T \bar{\mathbf H} \mathbf{e}_j + \mathbf{e}_j^T \bar{\mathbf H} \mathbf{e}_j \geq 2\mathbf{e}_i^T \bar{\mathbf H} \mathbf{e}_j + (1-\delta) (\|\mathbf e_i\|^2 + \|\mathbf e_j\|^2).
    	\end{equation}
    	Combining \eqref{eq:ec_ineq_1} and \eqref{eq:ec_ineq_2} yields that
    	\[
    		\mathbf{e}_i^T \bar{\mathbf H} \mathbf{e}_j \leq \delta(\|\mathbf e_i\|^2 + \|\mathbf e_j\|^2) \qquad \forall i \neq j.
    	\]
    	Therefore,
    	\begin{align*}
    		\hat{\mathbf e}^T \hat{\mathbf H} \hat{\mathbf e} = (\sum_{i=3}^l \mathbf{e}_i)^T \bar{\mathbf H} (\sum_{i=3}^l \mathbf{e}_i) &\leq (1+\delta)(\sum_{i=3}^l  \|\mathbf{e}_i\|^2) + 2 \delta (l-3)  (\sum_{i=3}^l \|\mathbf{e}_i\|^2) \\
    		&= (1+\delta)\|\mathbf{e}_c\|^2 + 2 \delta (l-3)  \|\mathbf{e}_c\|^2.
    	\end{align*}
    	The above inequality directly verifies the satisfaction of the first constraint. For the second constraint, consider an arbitrary vector $
    	\tilde{\mathbf{e}} \in \RR^{4r^2}$. Then,
    	\begin{align*}
    		\tilde{\mathbf{e}}^T \hat{\mathbf H} \tilde{\mathbf{e}} = \tilde{\mathbf{e}}^T \mathbf{P}^T \bar{\mathbf H} \mathbf{P} \tilde{\mathbf{e}} &= \tilde{\mathbf{e}}^T (P^T \otimes P^T) \bar{\mathbf H} (P \otimes P) \tilde{\mathbf{e}} \\
    		&= \vecc(P \mat(\tilde{\mathbf{e}}) P^T)^T \bar{\mathbf H} \vecc(P \mat(\tilde{\mathbf{e}}) P^T).
    	\end{align*}
    	By orthogonal projection, we know that $P \mat(\tilde{\mathbf{e}}) P^T \in \RR^{n \times n}$ has rank $2r$. Therefore, the following holds by the $\delta_{2r}$-RIP property of $\bar{\mathbf H}$:
    	\begin{equation}\label{eq:tilde_e_rip}
    		(1-\delta) \|P \mat(\tilde{\mathbf{e}}) P^T\|^2_F \leq \vecc(P \mat(\tilde{\mathbf{e}}) P^T)^T \bar{\mathbf H} \vecc(P \mat(\tilde{\mathbf{e}}) P^T) \leq (1+\delta) \|P \mat(\tilde{\mathbf{e}}) P^T\|^2_F
    	\end{equation}
    	and since
    	\begin{align*}
    		\|P \mat(\tilde{\mathbf{e}}) P^T\|^2_F &= \tr(P \mat(\tilde{\mathbf{e}})^T P^T P \mat(\tilde{\mathbf{e}}) P^T) \\
    		&= \tr(P \mat(\tilde{\mathbf{e}})^T \mat(\tilde{\mathbf{e}}) P^T) \\
    		&= \tr(P^T P \mat(\tilde{\mathbf{e}})^T \mat(\tilde{\mathbf{e}})) \\
    		&= \tr(\mat(\tilde{\mathbf{e}})^T \mat(\tilde{\mathbf{e}})) \\
    		&= \|\tilde{\mathbf{e}}\|^2_2,
    	\end{align*}
    	\eqref{eq:tilde_e_rip} automatically implies the satisfaction of the second constraint.
    \end{proof}

    \begin{proof}[Proof of Lemma \ref{lem:delta_ineq_2}] \renewcommand{\qedsymbol}{}
    	It suffices to show that for any feasible pair $(\delta, \hat{\mathbf H})$ of \eqref{eq:int_lmi}, we can construct a feasible solution
    	$(\delta,\mathbf H)$ to \eqref{eq:lmi} characterized as
    	\[
    		\delta = \delta, \qquad \mathbf H = \mathbf{P} \hat{\mathbf H} \mathbf{P}^T + (1-\delta)(\mathbf{I}_{n^2} - \mathbf{P}\mathbf{P}^T).
    	\]
    	To prove the lemma, it is enough to verify that the above pair $(\delta, \mathbf H)$ is feasible to \eqref{eq:lmi}. We first verify the second constraint. Given an arbitrary vector $\mathbf{e} \in \RR^{n^2}$, we have that
    	\[
    		\mathbf{e}^T \mathbf H \mathbf{e} = \mathbf{e}^T \mathbf{P} \hat{\mathbf H} \mathbf{P}^T \mathbf{e} + (1-\delta) \left[ \mathbf{e}^T \mathbf{e} - \mathbf{e}^T \mathbf{P}\mathbf{P}^T \mathbf{e} \right]
    	\]
    	and defining $\tilde{\mathbf{e}} \coloneqq \mathbf{P}^T \mathbf{e} \in \RR^{4r^2}$, we obtain:
    	\[
    		\mathbf{e}^T \mathbf{P} \hat{\mathbf H} \mathbf{P}^T \mathbf{e} + (1-\delta) \left[ \mathbf{e}^T \mathbf{e} - \mathbf{e}^T \mathbf{P}\mathbf{P}^T \mathbf{e} \right] \geq (1-\delta)\|\tilde{\mathbf{\mathbf{e}}}\|^2_2 + (1-\delta) [\|\mathbf{e}\|^2_2 - \|\tilde{\mathbf{\mathbf{e}}}\|^2_2] = (1-\delta)\|\mathbf{e}\|^2_2.
    	\]
    	Also, since $\|\tilde{\mathbf{e}}\|^2_2 \leq \|\mathbf{e}\|^2_2$ and $P$ is a projection matrix, one can write:
    	\[
    		(1+\delta)[\|\mathbf{e}\|^2_2 - \|\tilde{\mathbf{e}}\|^2_2] \geq (1-\delta)[\|\mathbf{e}\|^2_2 - \|\tilde{\mathbf{e}}\|^2_2],
    	\]
    	which further implies that
    	\[
    		\mathbf{e}^T \mathbf{P} \hat{\mathbf H} \mathbf{P}^T \mathbf{e} + (1-\delta) \left[ \mathbf{e}^T \mathbf{e} - \mathbf{e}^T \mathbf{P}\mathbf{P}^T \mathbf{e} \right] \leq (1+\delta)\|\tilde{\mathbf{e}}\|^2_2 +(1+\delta)[\|\mathbf{e}\|^2_2 - \|\tilde{\mathbf{e}}\|^2_2] = (1+\delta) \|\mathbf{e}\|^2_2.
    	\]
    	Combining the above equations, we recover the second constraint of \eqref{eq:lmi}:
    	\[
    		(1-\delta)\|\mathbf{e}\|^2_2\leq e^T \mathbf H \mathbf{e} \leq (1+\delta)\|\mathbf{e}\|^2_2.
    	\]
    	To study the first constraint, we have that
    	\begin{align*}
    		\mathbf \mathbf{e}^T \mathbf H \mathbf \mathbf{e} &= \hat{\mathbf e}^T \hat{\mathbf H} \hat{\mathbf e} + (1-\delta) \left[ \|\mathbf e\|^2_2 - \| \hat{\mathbf e}\|^2_2 \right]\\
    		 &\leq (1+\delta) \|\mathbf{e}_c\|^2_2 + 2(l-3)\delta \|\mathbf{e}_c\|^2_2 + (1-\delta) \|\mathbf{e}_c\|^2_2\\
    		  &= 2 \|\mathbf{e}_c\|^2_2 + 2(l-3)\delta \|\mathbf{e}_c\|^2_2.
    	\end{align*}
    	Note that $\|\mathbf e\|^2_2 - \| \hat{\mathbf e}\|^2_2 = \|\mathbf{e}_c\|^2_2$ due to
    	\[
    		\|\mathbf e\|^2_2 = \sum_{i=1}^n \lambda_i^2, \qquad \|\hat{\mathbf e}\|^2_2 = \sum_{i=1}^{2r} \lambda_i^2.
    	\]
    \end{proof}
The proof of Proposition \ref{prop:convex_lb} follows directly from combining Lemma \ref{lem:delta_ineq_1} and \ref{lem:delta_ineq_2}.

\end{proof}

\subsection{Proof of Lemma \ref{lem:delta_lb}}

\begin{proof}
    We aim to solve \eqref{eq:lmi} analytically to obtain a sufficient RIP bound for problem \eqref{eqn:ms-mc-sdp}. This amounts to deriving a closed-form expression for $\delta_{\text{lb}}(\mathbf e)$. We consider a simpler problem to solve \eqref{eq:lmi}:
    \begin{equation}\label{eqn:eta_lmi}
    \begin{aligned}
    \max_{\eta,\tilde{\mathbf H}} \quad & \eta \\
    \st \quad & \mathbf e^T \tilde{\mathbf H}\mathbf e \leq \frac{1-\eta}{2} c^2+\frac{1+\eta}{2} d^2 \\
    & \eta I_{n^2} \preceq \tilde{\mathbf H} \preceq I_{n^2}
    \end{aligned}
    \end{equation}
    with $c^2 = 2(l-3) \|\mathbf{e}_c\|^2_2$   and $d^2 = 2 \|\mathbf{e}_c\|^2_2$. Given any feasible solution $(\delta,\mathbf H)$ to \eqref{eq:lmi}, the tuple
    \[
    \left(\frac{1-\delta}{1+\delta},\frac{1}{1+\delta}\mathbf H\right)
    \]
    is a feasible solution to problem \eqref{eqn:eta_lmi}. Therefore, if we denote the optimal value of \eqref{eqn:eta_lmi} as $\eta(\mathbf e)$, then it holds that
    \begin{equation}\label{eqn:delta_eta}
    \eta(\mathbf e) \geq \frac{1-\delta_{\text{lb}}(\mathbf e)}{1+\delta_{\text{lb}}(\mathbf e)} \implies \delta_{\text{lb}}(\mathbf e) \geq \frac{1-\eta(\mathbf e)}{1+\eta(\mathbf e)}.
    \end{equation}
    We use the dual problem to solve for $\eta(\mathbf e)$:
    \begin{equation}\label{eqn:eta_lmi_dual}
    \begin{aligned}
    \min_{\mathbf{U}_1, \mathbf{U}_2,\gamma} \quad & \tr(\mathbf{U}_2) +\frac{\gamma}{2} (c^2+d^2)\\
    \st \quad & \tr(\mathbf{U}_1)+\frac{\gamma}{2}(c^2-d^2)=1 \\
    & \gamma \mathbf e \mathbf e^T = \mathbf{U}_1 - \mathbf{U}_2, \quad \mathbf{U}_1, \mathbf{U}_2 \succeq 0, \ \gamma \geq 0.
    \end{aligned}
    \end{equation}
    Since Slater's condition holds for the convex program \eqref{eqn:eta_lmi}, the optimal solution to \eqref{eqn:eta_lmi_dual} is equivalent to that of \eqref{eqn:eta_lmi}, which is $\eta(\mathbf e)$. Using a Lagrangian argument, $\eta(\mathbf e)$ can be solved as follows:
    \begin{align*}
    	\eta(\mathbf e) &=\max_{\beta \in \RR} \min_{\gamma \geq 0} \left\{  \beta(1-\frac{\gamma}{2}(c^2-d^2)) + \gamma \frac{c^2+d^2}{2} + \min_{\substack{\mathbf{U}_1 \succeq 0\\ \mathbf{U}_1 - \gamma \mathbf{e} \mathbf{e}^T \succeq 0}} \left[ \tr(\mathbf{U}_1 - \gamma \mathbf{e} \mathbf{e}^T ) - \beta \tr(\mathbf{U}_1 )\right] \right\} \\
    	&= \max_{\beta \leq 1} \min_{\gamma \geq 0} \left\{  \beta(1-\frac{\gamma}{2}(c^2-d^2)) + \gamma \frac{c^2+d^2}{2} + \min_{\substack{\mathbf{U}_1 \succeq 0\\ \mathbf{U}_1 - \gamma \mathbf{e} \mathbf{e}^T \succeq 0}} \left[ (1-\beta)\tr(\mathbf{U}_1) - \gamma \|\mathbf e\|^2_2 \right] \right\} \\
    	&= \max_{\beta \leq 1} \min_{\gamma \geq 0} \left\{  \beta(1-\frac{\gamma}{2}(c^2-d^2)) + \gamma \frac{c^2+d^2}{2} +  \gamma (1-\beta)\|\mathbf e\|^2_2 - \gamma \|\mathbf e\|^2_2 \right\} \\
    	&= \max_{\beta \leq 1}  \left\{ \beta + \min_{\gamma \geq 0} \left[ \gamma (\frac{c^2+d^2}{2} - \beta(\frac{c^2-d^2}{2}+\|\mathbf e\|^2_2))\right] \right\} \\
    	&=  \max_{\beta \leq 1}  \left\{ \beta: \frac{c^2+d^2}{2} - \beta(\frac{c^2-d^2}{2}+\|\mathbf e\|^2_2) \geq 0 \right\} \\
    	&= \min \left\{1, \frac{c^2+d^2}{2\|\mathbf e\|^2_2 + c^2 -d^2}\right\},
    \end{align*}
    where the first equality uses the Lagrangian argument by introducing the Lagrange multiplier $\beta$, and the second equality constraints $\beta \leq 1$ since $(1-\beta)\tr(\mathbf{U}_1)$ will be unbounded otherwise. The third equality results from the obvious choice of $\mathbf{U}_1 = \gamma \mathbf{e} \mathbf{e}^T$ given that $(1-\beta)$ is nonnegative. The fifth equality results from the choice of $\gamma = 0$ constrained to the requirement that its coefficient must be nonnegative.
    
    Substituting $\eta(\mathbf e) = 1$ into \eqref{eqn:delta_eta} results in the trivial lower bound $\delta_{\text{lb}}(\mathbf e)$ of 0, which means that the lower bound indeed will not be negative. Hence, we will focus on $\frac{c^2+d^2}{2\|\mathbf e\|^2_2 + c^2 -d^2}$ from now on. We know from \eqref{eqn:delta_eta} that in order to obtain a lower bound on $\delta_{\text{lb}}$, we need to derive an upper bound on $\eta(\mathbf e)$. Note that
    \begin{equation} \label{eqn:eta_e_analytic}
    	\frac{c^2+d^2}{2\|\mathbf e\|^2_2 + c^2 -d^2} = \frac{2(l-2)\|\mathbf{e}_c\|^2_2}{2\|\mathbf{e}\|^2_2+ 2(l-4)\|\mathbf{e}_c\|^2_2 } = \frac{(l-2) \|\mathbf{e}_c\|^2_2}{\|\mathbf{e}_{2r}\|^2_2 + (l-3) \|\mathbf{e}_c\|^2_2}.
    \end{equation}
    The last equality follows from the relations
    \[
    	\|\mathbf{e}\|^2_2 = \|\mathbf{e}_c +\mathbf{e}_{2r} \|^2_2 = \|\mathbf{e}_{2r}\|^2_2+ \|\mathbf{e}_{c}\|^2_2.
    \]
    If we fix $\|\mathbf{e}_{2r}\|^2_2$, then we can maximize \eqref{eqn:eta_e_analytic} with respect to $ \|\mathbf{e}_{c}\|^2_2$ first. In this case, taking the derivative of \eqref{eqn:eta_e_analytic} yields that
    \[
    	\frac{\partial}{\partial \|\mathbf e_c\|_2} \left( \frac{(l-2) \|\mathbf{e}_c\|^2_2}{\|\mathbf{e}_{2r}\|^2_2 + (l-3) \|\mathbf{e}_c\|^2_2} \right) = 2\frac{(l-2)\|\mathbf e_c\|_2 \|\mathbf{e}_{2r}\|^2_2}{(\|\mathbf{e}_{2r}\|^2_2 + (l-3) \|\mathbf{e}_c\|^2_2)^2} \geq 0.
    \]
    Therefore, \eqref{eqn:eta_e_analytic} is maximized when $\|\mathbf e_c\|_2^2$ is set to be as large as possible. Before we derive the maximum value of $\|\mathbf e_c\|_2^2$ in terms of $\|\mathbf e_1\|_2^2$ and $\|\mathbf e_2\|_2^2$, we introduce one key lemma.
    \begin{lemma}\label{lem:approx_lr}
	    Consider two PSD matrices $\mathbf{M}$ and $ \mathbf{M}^*$ such that $\tr(\mathbf{M}) \leq \tr(\mathbf{M}^*)$ and $\rk(\mathbf{M}^*) = r$. Then,
	    \begin{equation}\label{eqn:eigva_rel}
	    \sigma_{(1)}(\mathbf{M}^* - \mathbf{M}) + \cdots + \sigma_{(r)}(\mathbf{M}^* - \mathbf{M}) \geq \sigma_{(r+1)}(\mathbf{M}^* - \mathbf{M}) + \cdots \sigma_{(n)}(\mathbf{M}^* - \mathbf{M}),
        \end{equation}
        where $\sigma_{i}$ denote the $i$-th largest singular value  of the matrix $\mathbf{M}^* - \mathbf{M}$.
    \end{lemma}
    \begin{proof}
    	For each matrix $\mathbf{A}$, we denote  the $i^{\text{th}}$ eigenvalue as $\lambda_{(i)}(\cdot)$, meaning that
    	\[
    		\lambda_{(1)}(\mathbf{A}) \geq \lambda_{(2)}(\mathbf{A}) \geq \dots \geq \lambda_{(n)}(\mathbf{A}).
    	\]
    	By Weyl's inequality, we know that
    	\[
    		\lambda_{(i+j-1)} (\mathbf{M}^* - \mathbf{M}) \leq \lambda_{(i)}(\mathbf{M}^*) + \lambda_{(j)}(-\mathbf{M}).
    	\]
    	Hence,
    	\[
    		\lambda_{(r+1)} (\mathbf{M}^*-\mathbf{M}) \leq \lambda_{(r+1)}(\mathbf{M}^*) + \lambda_{(1)}(-\mathbf{M}) \leq 0
    	\]
    	since $\mathbf{M}^*$ is of rank-$r$ and $\mathbf{M} \succeq 0$. Therefore, we know that $\mathbf{M}^*-\mathbf{M}$ has at most $r$ positive eigenvalues. Also, since $\tr(\mathbf{M}^*-\mathbf{M}) \geq 0$, it holds that
    	\[
    		 \lambda_{(1)}(\mathbf{M}^*-\mathbf{M}) + \dots + \lambda_{(r)}(\mathbf{M}^*-\mathbf{M}) \geq - \lambda_{(r+1)}(\mathbf{M}^*-\mathbf{M}) - \dots - \lambda_{(n)}(\mathbf{M}^*-\mathbf{M}),
    	\]
    	which implies that 
    	\begin{equation}\label{eq:weyl_corollary}
    		|\lambda_{(1)}(\mathbf{M}^*-\mathbf{M})| + \dots + |\lambda_{(r)}(\mathbf{M}^*-\mathbf{M})| \geq |\lambda_{(r+1)}(\mathbf{M}^*-\mathbf{M})| + \dots + |\lambda_{(n)}(\mathbf{M}^*-\mathbf{M})|
    	\end{equation}
    	since $\lambda_{(k)}(\mathbf{M}^*-\mathbf{M}) \leq 0$ for all $k > r$. According to the definition, we have
    	\begin{equation}\label{eq:absolute_order}
    		\begin{aligned}
    			&|\lambda_1(\mathbf{M}^*-\mathbf{M})| + \dots + |\lambda_r(\mathbf{M}^*-\mathbf{M})| \geq |\lambda_{(1)}(\mathbf{M}^*-\mathbf{M})| + \dots + |\lambda_{(r)}(\mathbf{M}^*-\mathbf{M})|, \\
    			& |\lambda_{(r+1)}(\mathbf{M}^*-\mathbf{M})| + \dots + |\lambda_{(n)}(\mathbf{M}^*-\mathbf{M})| \geq |\lambda_{r+1}(\mathbf{M}^*-\mathbf{M})| + \dots |\lambda_n(\mathbf{M}^*-\mathbf{M})|
    		\end{aligned}
    	\end{equation}
    	since $\lambda_1(\mathbf{M}^*-\mathbf{M}), \dots, \lambda_n(\mathbf{M}^*-\mathbf{M})$ are ordered with respect to their absolute values. As per the main text, we abbreviate $\lambda_i(\mathbf{M}^*-\mathbf{M})$ as $\lambda_i$ for the sake of brevity for any $i \in [n]$. Combining \eqref{eq:weyl_corollary}  with  \eqref{eq:absolute_order} proves the original lemma.
    \end{proof}
    
    Denote $S_1 \coloneqq \sum_{i=1}^r |\lambda_i|$ and $ S_2 \coloneqq \sum_{i=r+1}^n |\lambda_i|$. Given $S_2$, since $|\lambda_i| \leq |\lambda_r|$ as long as $ i > r$, we know that $\sum_{i=r+1}^n \lambda_i^2$ is maximized when every absolute value is chosen to be as large as possible, namely
    \[
    	|\lambda_{r+1}|, \dots, |\lambda_{r+ \lfloor S_2/|\lambda_r| \rfloor}| = |\lambda_r|, |\lambda_{r+ \lceil S_2/|\lambda_r| \rceil}| = S_2 - \lfloor S_2/|\lambda_r| \rfloor  |\lambda_r| \coloneqq \tilde \lambda \leq |\lambda_r|.
    \]
    Therefore,
    \[
    	\sum_{i=r+1}^n \lambda_i^2 \leq \lfloor S_2/|\lambda_r| \rfloor \lambda_r^2 + \tilde \lambda^2 \leq \lfloor S_2/|\lambda_r| \rfloor \lambda_r^2 + \frac{\tilde \lambda}{|\lambda_r|}  \lambda_r^2 = \frac{S_2}{|\lambda_r|} \lambda_r^2.
    \]
    As a result,
    \begin{align*}
    	\frac{S_2}{|\lambda_r|} \lambda_r^2 = S_2 |\lambda_r| \leq S_1 |\lambda_r| \leq S_1 \frac{S_1}{r} \leq \frac{S^2_1}{r} \leq \frac{r \|\mathbf e_1\|^2_2}{r} =\|\mathbf e_1\|^2_2,
    \end{align*}
    where the last inequality follows from Cauchy-Schwartz. Combining the above 2 inequalities, we obtain
    \begin{equation}\label{eqn:e_normsq_rela}
    	\sum_{i=r+1}^n \lambda_i^2 \leq \|\mathbf e_1\|^2_2. 
    \end{equation}
    Furthermore, since $|\lambda_{r+1}| \geq \dots \geq |\lambda_{n}|$, one can write:
    \[
    	\|\mathbf e_c\|^2_2 = \sum_{i=2r+1}^n \lambda_i^2 \leq \frac{n-2r}{n-r} \sum_{i=r+1}^n \lambda_i^2  
    \]
    with equality holding if and only if $|\lambda_{r+1}| = \dots = |\lambda_{n}|$.  Combined with \eqref{eqn:e_normsq_rela}, we obtain
    \begin{equation}\label{eqn:ec_bound}
    	\|\mathbf e_c\|^2_2 \leq \frac{n-2r}{n-r} \|\mathbf e_1\|^2_2.
    \end{equation}
    Consequently,
    \begin{equation} \label{eqn:e2_bound}
    	\|\mathbf e_c\|^2_2 \leq   \frac{n-2r}{n-r} (\|\mathbf e_2\|^2_2 + \|\mathbf e_c\|^2_2) \implies \|\mathbf e_2\|^2_2 \geq \frac{r}{n-2r} \|\mathbf e_c\|^2_2.
    \end{equation}
    It results from \eqref{eqn:ec_bound} and \eqref{eqn:e2_bound} that
    \begin{align*}
    	\max_{\mathbf e: \tr(M) \leq \tr(M^*)} \eta(\mathbf e) &= \max_{\mathbf e: \tr(M) \leq \tr(M^*)} \frac{(l-2) \|\mathbf{e}_c\|^2_2}{\|\mathbf{e}_{2r}\|^2_2 + (l-3) \|\mathbf{e}_c\|^2_2}\\
    	 &\leq \frac{(l-2) \|\mathbf{e}_c\|^2_2}{\|\mathbf{e}_{1}\|^2_2 +\frac{r}{n-2r} \|\mathbf e_c\|^2_2 + (l-3) \|\mathbf{e}_c\|^2_2} \\
    	&\leq \frac{(l-2)\frac{n-2r}{n-r}}{1+ (\frac{r}{n-2r}+l-3)\frac{n-2r}{n-r}}\\
    	&= \frac{(l-2)(n-2r)}{n+(n-2r)(l-3)}.
    \end{align*}
    Thus,
    \begin{align*}
    	\delta_{\text{lb}} \geq \max_{\mathbf e: \tr(M) \leq \tr(M^*)}  \frac{1-\delta(\mathbf e)}{1+\delta(\mathbf e)} = \frac{2r}{n+(n-2r)(2l-5)}.
    \end{align*}

\end{proof}

\subsection{Proof of Theorem \ref{thm: bm-solves-chain}}

\begin{proof}
    The matrix completion problem with the least-squares objective function can be written as 
\begin{equation} \label{eqn:chain}
    \min_{\mathbf{x} \in \mathbb{R}^n} f(\mathbf{x}) = (x_1^2 - (x_1^*)^2)^2 + 2 \sum_{i=1}^{n-1} (x_i x_{i+1} - x^*_i x^*_{i+1} )^2.
\end{equation}
We assume that every element of the vector $\mathbf{x}^*$ is nonzero, i.e. $x^*_i \not = 0, \forall i \in [n]$. In this case, the problem \eqref{eqn:chain} does not have any spurious solutions. In order to prove the inexistence of a non-global second-order critical solution, we need to investigate the second-order necessary optimality conditions. The gradient and Hessian of the above objective can be written as
\begin{align*}
    [\nabla f(x)]_i & = \begin{cases}
    (x_1^2 - (x^*_1)^2)x_1 + (x_1x_2 - x^*_1 x^*_2)x_2 , & \text{if }i =1\\
    (x_{i-1} x_{i} - x^*_{i-1} x^*_{i} )x_{i-1} + (x_i x_{i+1} - x^*_i x^*_{i+1} )x_{i+1}  , & \text{if }i \not \in  \{1,n\} \\
    (x_{n-1}x_n - x^*_{n-1} x^*_n)x_{n-1}, & \text{if }i =n
    \end{cases}, \\
    [\nabla^2 f(x)]_{i,j} & =  \begin{cases}
    3x_1^2 - (x^*_1)^2 + x_2^2 , & \text{if }i,j =1\\
    x_{i-1}^2 + x_{i+1}^2  , & \text{if }i=j,i,j \not \in  \{1,n\} \\
    x_{n-1}^2, & \text{if }i,j =n \\ 
    2x_ix_j - x^*_i x^*_j , & \text{if } |i-j| = 1\\
    0  , & \text{otherwise }
    \end{cases}
\end{align*}
Note that for the point $\mathbf{\hat{x}}$ to be a second-order stationary point, we require $\nabla f(\mathbf{\hat{x}}) = 0$ and $\nabla^2 f(\mathbf{\hat{x}}) \succeq 0$. $[\nabla f(\mathbf{\hat{x}})]_n = 0$ implies either $\hat{x}_{n-1} \hat{x}_n = x^*_{n-1} x^*_n $ or $\hat{x}_{n-1} = 0$. The latter results in $\nabla^2 f(\mathbf{\hat{x}}) \not \succeq 0$ since the following principal minor is not positive semidefinite: 
\begin{align*}
    [\nabla^2 f(\mathbf{\hat{x}})]_{n-1:n,n-1:n} & = \begin{bmatrix} (\hat{x}_{n-2})^2 + (\hat{x}_n)^2 & 2 \hat{x}_{n-1} \hat{x}_n - x^*_{n-1} x^*_n \\ 2 \hat{x}_{n-1} \hat{x}_n - x^*_{n-1} x^*_n & (\hat{x}_{n-1})^2\end{bmatrix} \\
    & = \begin{bmatrix} (\hat{x}_{n-2})^2 + (\hat{x}_n)^2 &  - x^*_{n-1} x^*_n \\  - x^*_{n-1} x^*_n & 0\end{bmatrix} \not \succeq 0.
\end{align*}
Hence, $\hat{x}_{n-1} \hat{x}^*_n = x^*_{n-1} x^*_n $ is the only feasible option for $\mathbf{\hat{x}}$ to be a second-order critical point. Next, we show that no second-order critical point $\mathbf{\hat{x}}$ can have a zero entry, i.e. $\hat{x}_i \not = 0 $ for all $i \in [n]$.

The first goal is to identify the structure of the stationary points with at least one zero entry. If $\hat{x}_1 = 0$, $[\nabla f(\mathbf{\hat{x}})]_1$ implies $\hat{x}_2 = 0$. Then, $[\nabla f(\mathbf{\hat{x}})]_2$ gives $\hat{x}_3 = 0$.  Continuing iteratively gives that the only possible stationary point with $\hat{x}_1 = 0$ is $\mathbf{\hat{x}} = 0$, which cannot be a second-order critical point. Similarly, if $\hat{x}_n = 0$, then $[\nabla f(\mathbf{\hat{x}})]_n$ implies $\hat{x}_{n-1} = 0$. Then, $[\nabla f(\mathbf{\hat{x}})]_{n-1}$ gives $\hat{x}_{n-2} = 0$. Hence, the only stationary point with $\hat{x}_n = 0$ is $\mathbf{\hat{x}} = 0$. Consequently, $\hat{x}_1 \not = 0$ and $\hat{x}_n \not = 0$ for a second-order stationary point.

Let $k$ be defined as an index with the property that $\hat{x}_1, \hat{x}_2, \dots, \hat{x}_{k-1} \not = 0$ and $\hat{x}_k = 0$. Then, $[\nabla f(\mathbf{\hat{x}})]_{k-2}$ yields $\hat{x}_{k-2} \hat{x}_{k-1} = x^*_{k-2} x^*_{k-1}$. Continuing iteratively backwards on $[\nabla f(\hat{\mathbf{x}}))]_{i}$ for $1 \leq i \leq k-2$ gives the following conditions on second-order critical points, which are correct values for the corresponding edges:
\begin{align*}
    \hat{x}_{k-2} \hat{x}_{k-1} & = x^*_{k-2} x^*_{k-1}, \\
    \hat{x}_{k-3} \hat{x}_{k-2} & = x^*_{k-3} x^*_{k-2}, \\
    & \vdots \\
    \hat{x}_{1} \hat{x}_{2} & = x^*_{1} x^*_{2} ,\\
    (\hat{x}_{1})^2 & = (x^*_1)^2.
\end{align*}
We can focus on entries of the stationary point $\mathbf{\hat{x}}$ corresponding to $\hat{x}_i$ for $i > k$. We show that $\hat{x}_{k+1} \not = 0$ and $\hat{x}_{k+2} \not = 0$ for every second-order critical point. If $\hat{x}_{k+1} = 0$, then $[\nabla f(\hat{\mathbf{x}})]_k = 0$ gives $x^*_{k-1} x^*_k = 0$, which contradicts the assumption on the ground truth matrix. Since $\hat{x}_{n-1}, \hat{x}_n \not = 0$, we have $k \leq n-2$ for each second-order critical point. If $\hat{x}_{k+2} = 0$, then the $2 \times 2$ submatrix of the Hessian will have the form
\begin{align*}
    [\nabla^2 f(\mathbf{\hat{x}})]_{k:k+1, k:k+1} & = \begin{bmatrix}
        (\hat{x}_{k-1})^2 + (\hat{x}_{k+1})^2 & 2 \hat{x}_k \hat{x}_{k+1} - x^*_k x^*_{k+1} \\ 2 \hat{x}_k \hat{x}_{k+1} - x^*_k x^*_{k+1} & (\hat{x}_{k})^2 + (\hat{x}_{k+2})^2
    \end{bmatrix} \\
    & = \begin{bmatrix}
        (\hat{x}_{k-1})^2 + (\hat{x}_{k+1})^2 &  - x^*_k x^*_{k+1} \\  - x^*_k x^*_{k+1} & 0
    \end{bmatrix} \not \succeq 0.
\end{align*}
Let $m$ be defined as an index with the property that $\hat{x}_{k+1}, \hat{x}_{k+2}, \dots, \hat{x}_{m} \not = 0$ and either $\hat{x}_{m+1} = 0$ or $m=n$. By the previous arguments and the definition of $m$, we have the condition $k+3 \leq m \leq n$. We are not interested in the entries after the $m$-th entry because we can show that a stationary point with this structure cannot be a second-order critical point because $[\nabla^2 f(\mathbf{\hat{x}})]_{1:m, 1:m} \not \succeq 0$. By using the first-order partial derivatives and algebra, we obtain the following equations for the stationary point $\mathbf{\hat{x}}$: 
    \begin{align*}
        & [\nabla f(\mathbf{\hat{x}})]_k = 0  \xrightarrow{} \hat{x}_{k+1} = - \Big(\frac{x^*_{k-1}}{x^*_{k+1}}\Big)^2 x^*_{k+1} = - \alpha x^*_{k+1}, \\
       &  [\nabla f(\mathbf{\hat{x}})]_{k+1} = 0  \xrightarrow{} \hat{x}_{k+2} = - \Big(\frac{x^*_{k-1}}{x^*_{k+1}}\Big)^{-2}x^*_{k+2} = - \alpha^{-1} x^*_{k+2}, \\
        & \vdots \\
        & [\nabla f(\mathbf{\hat{x}})]_{m-1} = 0  \xrightarrow{} \hat{x}_{m} = - \Big(\frac{x^*_{k-1}}{x^*_{k+1}}\Big)^{2(-1)^{(m-k+1)}}x^*_{m} = - \alpha^{({(-1)}^{m-k+1})} x^*_{m},
    \end{align*}
    where $\alpha = \Big(\frac{x^*_{k-1}}{x^*_{k+1}}\Big)^2 $. As a result, the possible candidates for second-order critical points with at least one zero entry have the form:
    \begin{align*}
        \mathbf{\hat{x}} & = [ x^*_1, x^*_2, \dots, x^*_{k-1}, 0, -\alpha x^*_{k+1}, -\alpha^{-1} x^*_{k-2}, \dots, - \alpha^{({(-1)}^{m-k+1})} x^*_{m}, \dots] 
    \end{align*}
    or the form
    \begin{align*}
        \mathbf{\hat{x}} & = [ -x^*_1, -x^*_2, \dots, -x^*_{k-1}, 0, \alpha x^*_{k+1}, \alpha^{-1} x^*_{k-2}, \dots,  \alpha^{({(-1)}^{m-k+1})} x^*_{m}, \dots].
    \end{align*}
    These points correspond to the same matrix completion and they lead to the same Hessian matrix. As mentioned before, we focus on the $m \times m $ Hessian submatrix $[\nabla^2 f(\hat{\mathbf{x}})]_{1:m, 1:m}$:
    \begin{align*}
        [\nabla^2 f(\mathbf{\hat{x}})]_{1:m, 1:m} = \begin{pmatrix}
          \mathbf{E}
          & \rvline & \mathbf{D} \\
        \hline
          \mathbf{D}^T & \rvline &
          \begin{matrix}
          \mathbf{A} & \mathbf{B} \\
          \mathbf{B}^T & \mathbf{C}
          \end{matrix}
        \end{pmatrix},
    \end{align*}
    where $\mathbf{A} \in \mathbb{R}^{3 \times 3}, \mathbf{B} \in \mathbb{R}^{3 \times (m-k-1)}, \mathbf{C} \in \mathbb{R}^{(m-k-1) \times (m-k-1)}, \mathbf{D} \in \mathbb{R}^{(k-2) \times (m-k+2)}$ and $\mathbf{E} \in \mathbb{R}^{(k-2) \times (k-2)}$. The submatrices can be written as
    \begin{align*}
        \mathbf{A} & =
        \begin{bmatrix}
            (x^*_{k-2})^2 & -x^*_{k-1} x^*_{k} & 0 \\
            - x^*_{k-1} x^*_{k} & (x^*_{k-1})^2 + \frac{(x^*_{k-1})^4}{(x^*_{k+1})^2} & - x^*_{k} x^*_{k+1} \\
            0 & - x^*_{k}x^*_{k+1} & \frac{(x^*_{k+1})^4 (x^*_{k+2})^2}{(x^*_{k-1})^4}
        \end{bmatrix}, \\
        \mathbf{B} & = \begin{bmatrix}
            0 & 0 & \dots & 0 \\
            0 & 0 & \dots & 0 \\
            x^*_{k+1} x^*_{k+2} & 0 & \dots & 0 
        \end{bmatrix}, \\
        \mathbf{C} & = \begin{bmatrix}
            \alpha^2 (x^*_{k+1})^2 + \alpha^2 (x^*_{k+3})^2 & x^*_{k+2}x^*_{k+3} & 0 & \dots & 0 \\
            x^*_{k+2}x^*_{k+3} & \alpha^{-2} (x^*_{k+2})^2 + \alpha^{-2} (x^*_{k+4})^2 & x^*_{k+3}x^*_{k+4} & \dots & 0 \\
            0 & x^*_{k+3}x^*_{k+4} &  \alpha^2 (x^*_{k+3})^2 + \alpha^2 (x^*_{k+5})^2 & \dots & 0 \\
            \vdots & \vdots & \vdots & \ddots & \vdots \\
            0 & 0 & 0 & \dots & x^*_{m-1}x^*_{m} \\
            0 & 0 & 0 & \dots & (\alpha^2)^{({(-1)}^{m-k})} (x^*_{m-1})^2
        \end{bmatrix}, \\
        \mathbf{D} & = \begin{bmatrix}
            0 & 0 & \dots & 0 \\
            \vdots & \vdots &  \ddots & \vdots \\
            x^*_{k-2} x^*_{k-1} & 0 & \dots & 0
        \end{bmatrix}, \\
        \mathbf{E} & = [\nabla^2 f(\mathbf{x}^*)]_{1:k-2, 1:k-2}.
    \end{align*}
    We investigate three different cases: $\mathbf{A} \not \succeq 0$, $\mathbf{A} \succeq 0$ with at least one of the eigenvalues being equal to $0$, and $\mathbf{A} \succ 0$. If $\mathbf{A} \not \succeq 0$, then $ [\nabla^2 f(\mathbf{\hat{x}})]_{1:m, 1:m} \not \succeq 0$ because $\mathbf{A}$ is a principal minor of the Hessian. Therefore, $\mathbf{\hat{x}}$ cannot be a second-order critical point. If $\mathbf{A} \succeq 0$ with an eigenvalue equal to $0$, we consider the following $4 \times 4$ principal minors of the Hessian:
     \begin{align*}
        \begin{pmatrix}
           (x^*_{k-3})^2 + (x^*_{k-1})^2
          & \rvline & \begin{matrix}
          x^*_{k-2} x^*_{k-1} & 0 & 0
          \end{matrix}  \\
        \hline
          \begin{matrix}
          x^*_{k-2} x^*_{k-1} \\ 0 \\ 0
          \end{matrix}  & \rvline &
          \mathbf{A}
        \end{pmatrix}, \qquad  \begin{pmatrix}
           \mathbf{A}
          & \rvline & \begin{matrix}
          0 \\ 0  \\ x^*_{k+1} x^*_{k+2}
          \end{matrix}  \\
        \hline
          \begin{matrix}
          0 & 0 & x^*_{k+1} x^*_{k+2}
          \end{matrix}  & \rvline &
          \alpha^2 (x^*_{k+1})^2 + \alpha^2 (x^*_{k+3})^2
        \end{pmatrix}.
    \end{align*}
    The submatrices are not positive semidefinite by Schur complement if
    \begin{align*}
        \mathbf{A} - \frac{1}{(x^*_{k-3})^2 + (x^*_{k-1})^2} \begin{bmatrix}
            (x^*_{k-2})^2 (x^*_{k-1})^2 & 0 & 0 \\ 0 & 0 & 0 \\ 0 & 0 & 0
        \end{bmatrix} \not \succeq 0
    \end{align*}
    and
    \begin{align*}
        \mathbf{A} - \frac{1}{\alpha^2 (x^*_{k+1})^2 + \alpha^2 (x^*_{k+3})^2} \begin{bmatrix}
            0 & 0 & 0 \\ 0 & 0 & 0 \\ 0 & 0 & (x^*_{k+1})^2 (x^*_{k+2})^2
        \end{bmatrix} \not \succeq 0,
    \end{align*}
    respectively. If $\mathbf{A}\mathbf{v} =0$, then $v_1$ and $v_3$ cannot be equal to $0$ at the same time because $(x^*_{k-1})^2 + \frac{(x^*_{k-1})^4}{(x^*_{k+1})^2} > 0$. If $v_1 \not = 0$, then 
    \begin{align*}
        \mathbf{v}^T\mathbf{A}\mathbf{v} - \mathbf{v}^T\frac{1}{(x^*_{k-3})^2 + (x^*_{k-1})^2} \begin{bmatrix}
            (x^*_{k-2})^2 (x^*_{k-1})^2 & 0 & 0 \\ 0 & 0 & 0 \\ 0 & 0 & 0
        \end{bmatrix}\mathbf{v} = 0 - \frac{v_1^2 (x^*_{k-2})^2 (x^*_{k-1})^2}{(x^*_{k-3})^2 + (x^*_{k-1})^2} < 0.
    \end{align*}
    Hence, the quadratic form of the Hessian has a descent direction. If $v_1 = 0$, then $v_3 \not =0$ and we can use the second submatrix to show that the submatrix is not positive semidefinite. Thus, $\mathbf{\hat{x}}$ cannot be a second-order critical point. As a result, if $\mathbf{\hat{x}}$ is a second-order critical point, then $\mathbf{A} \succ 0$. In that case, we can show that the submatrix 
    \begin{align*}
        \begin{bmatrix}
            \mathbf{A} & \mathbf{B} \\
            \mathbf{B}^T & \mathbf{C}
        \end{bmatrix} 
    \end{align*}
    cannot be positive semidefinite. This requires proving that $\mathbf{C} - \mathbf{B}^T \mathbf{A}^{-1} \mathbf{B} \not \succeq 0$. Note that
    \begin{align*}
        \mathbf{B}^T \mathbf{A}^{-1} \mathbf{B} = \begin{bmatrix}
            \mathbf{A}^{-1}_{3,3} (x^*_{k+1})^2 (x^*_{k+2})^2 & 0 & \dots & 0 \\
            0 & 0 & \dots & 0 \\
            \vdots & \vdots & \ddots & \vdots \\
            0 & 0 & \dots & 0
        \end{bmatrix}.
    \end{align*}
    We can calculate $\mathbf{A}_{3,3}^{-1}$ by using the cofactors of $\mathbf{A}$ as
    \begin{align*}
        \mathbf{A}_{3,3}^{-1} = \frac{1}{\text{det}(\mathbf{A})} \Big( (x^*_{k-2})^2 (x^*_{k-1})^2 + (x^*_{k-2})^2 \frac{(x^*_{k-1})^4}{(x^*_{k+1})^2} - (x^*_{k-1})^2 (x^*_{k})^2 \Big),
    \end{align*}
    where 
    \begin{align*}
        \text{det}(\mathbf{A}) = (x^*_{k-2})^2 \Big( \frac{(x^*_{k+1})^4}{(x^*_{k-1})^2} (x^*_{k+2})^2 + (x^*_{k+1})^2 (x^*_{k+2})^2 - (x^*_{k})^2 (x^*_{k+1})^2 \Big) - \frac{(x^*_{k})^2 (x^*_{k+1})^4 (x^*_{k+2})^2}{(x^*_{k-1})^2}.
    \end{align*}
    An algebraic manipulation shows that $\frac{d \mathbf{A}_{3,3}^{-1}}{d (x^*_k)^2} > 0$. Thus, the value of $\mathbf{A}_{3,3}^{-1}$ is minimized when $(x^*_k)^2 \xrightarrow{} 0$. Whenever $(x^*_k)^2 = 0$, $\mathbf{A}^{-1}_{3,3} (x^*_{k+1})^2 (x^*_{k+2})^2$ is equal to $\alpha^2 (x^*_{k+1})^2$. As a result, $\mathbf{A}^{-1}_{3,3} (x^*_{k+1})^2 (x^*_{k+2})^2 > \alpha^2 (x^*_{k+1})^2$ by the non-zero assumption on $x^*_k$. We can write the matrix $\mathbf{C} - \mathbf{B}^T \mathbf{A}^{-1} \mathbf{B}$ as the summation of $m-k-1$ matrices as
    \begin{align*}
        \mathbf{C} - \mathbf{B}^T \mathbf{A}^{-1} \mathbf{B} = & \begin{bmatrix}
            \alpha^2 (x^*_{k+1})^2 - \mathbf{A}_{3,3}^{-1} (x^*_{k+1})^2 (x^*_{k+2})^2 & 0 & \dots & 0 & 0 \\
            0 & 0 & \dots & 0 & 0 \\
            \vdots & \vdots & \ddots & \vdots & \vdots  \\
            0 & 0 & \dots & 0 & 0 \\
            0 & 0 & \dots & 0 & 0
        \end{bmatrix} + \\
        & \begin{bmatrix}
            \alpha^2 (x^*_{k+3})^2 & x^*_{k+2} x^*_{k+3} & \dots & 0 & 0 \\
            x^*_{k+2} x^*_{k+3} & \alpha^{-2} (x^*_{k+2})^2 & \dots & 0 & 0 \\
            \vdots & \vdots & \ddots & \vdots & \vdots  \\
            0 & 0 & \dots & 0 & 0 \\
            0 & 0 & \dots & 0 & 0
        \end{bmatrix} + \cdots + \\
        & \begin{bmatrix}
            0 & 0 & \dots & 0 & 0 \\
            0 & 0 & \dots & 0 & 0 \\
            \vdots & \vdots & \ddots & \vdots & \vdots  \\
            0 & 0 & \dots & (\alpha^2)^{(-1)^{m-k+1}} (x^*_m)^2 & x^*_{m-1} x^*_m \\
            0 & 0 & \dots & x^*_{m-1} x^*_m & (\alpha^2)^{(-1)^{m-k}} (x^*_{m-1})^2 
        \end{bmatrix}.
    \end{align*}
    Consider the vector $\mathbf{v}$ defined as $v_{2t+1} = \frac{x^*_{k + 2t + 2}}{x^*_{k+2}}$ and $v_{2t} = \alpha^2 \frac{x^*_{k+2t+1}}{x^*_{k+2}}$ for $t = 0, \dots, \lfloor (m-k-1)/2 \rfloor$. Then,
    \begin{align*}
        \mathbf{v}^T (\mathbf{C} - \mathbf{B}^T \mathbf{A}^{-1} \mathbf{B}) \mathbf{v} = \alpha^2 (x^*_{k+1})^2 - \mathbf{A}_{3,3}^{-1} (x^*_{k+1})^2 (x^*_{k+2})^2 + 0 + \dots + 0 < 0.
    \end{align*}
    As a result, $ \mathbf{C} - \mathbf{B}^T \mathbf{A}^{-1} \mathbf{B}$ is not positive definite. Hence, the Hessian is not positive definite either. Consequently, the problem cannot have a second-order critical point $\hat{x}$ with $\hat{x}_i = 0$ for some $i \in [n]$. Then, the only possible second-order-critical points must satisfy $(\hat{x}_1)^2 = (x^*_1)^2$ and $\hat{x}_{i} \hat{x}_{i+1} = x^*_i x^*_{i+1}, i = 1, \dots, n-1$, which correspond to the valid factors of the ground truth completion of the matrix $\mathbf{M}^*$.
    
\end{proof}

\subsection{Proof of Theorem \ref{thm: sdp-fails-chain}}

\begin{proof}
     To show that the SDP formulated as \eqref{eqn:ms-mc-sdp} fails to solve the problem, consider two indices $j,k$ such that $x^*_k > x^*_j$ and $j,k > 2$. Then, we construct a feasible solution $\mathbf{\hat{M}}$ that is strictly better than the ground truth solution, which shows the failure of SDP. Let $\mathbf{\hat{M}} = \mathbf{y}\mathbf{y}^T + \mathbf{zz}^T$, sum of two rank-1 matrices, where
\[ y_i = x^*_i, \forall i \in [n] \backslash \{ k \}, \qquad y_k = x^*_j \]
and 
\[ z_i = 0, \forall i \in [n] \backslash \{ j,k \}, \qquad z_j = z_k = (|x^*_j x^*_k| - (x^*_j)^2)^{1/2}. \]
Since the sum of PSD matrices is PSD and $\mathbf{\hat{M}}$ satisfies the observed entries, $\mathbf{\hat{M}}$ is a feasible solution. Moreover, the objective value corresponding to $\mathbf{\hat{M}}$ is $\sum_{i \not = j, i \not = k} (x^*_i)^2 + 2|x^*_j x^*_k|$. Hence,  by the assumption, the feasible solution $\mathbf{\hat{M}}$ is strictly better than the ground truth solution $\mathbf{M}^*$. Thus, SDP fails to recover the true solution.
\end{proof}

\subsection{Proof of Theorem \ref{thm: bm-solves-cycle}}

\begin{proof}

Note that the nodes are numbered from $0$ to $2k$ as opposed to earlier examples. Hence, the matrix completion problem with the least-squares objective function can be written as
\begin{equation} \label{eqn:cycle}
    \min_{\mathbf{x} \in \mathbb{R}^n} f(\mathbf{x}) = 2 \sum_{i=0}^{2k} (x_i x_{i+1} - x^*_i x^*_{i+1} ) ^2 + 2 (x_0 x_{2k} - x^*_0 x^*_{2k})^2.
\end{equation}
The gradient and Hessian of the above objective can be written as
\begin{align*}
    [\nabla f(\mathbf{x})]_i & = \begin{cases}
    (x_0 x_{2k} - x^*_0 x^*_{2k})x_{2k} + (x_0 x_1 - x^*_0 x^*_1)x_1 , & \text{if }i =0\\
    (x_{i-1} x_{i} - x^*_{i-1} x^*_{i} )x_{i-1} + (x_i x_{i+1} - x^*_i x^*_{i+1} )x_{i+1}  , & \text{if }i \not \in  \{0,2k\} \\
    (x_{2k-1} x_{2k} - x^*_{2k-1} x^*_{2k}) x_{2k-1} + (x_0 x_{2k} - x^*_0 x^*_{2k} )x_0, & \text{if }i =2k
    \end{cases}, \\
    [\nabla^2 f(x)]_{i,j} & =  \begin{cases}
    x_1^2 + x_{2k}^2 , & \text{if }i,j =0\\
    x_{i-1}^2 + x_{i+1}^2  , & \text{if }i=j, \: i,j \not \in  \{0,2k\} \\
    x_0^2 + x_{2k-1}^2 , & \text{if }i,j =2k \\ 
    2x_ix_j - x^*_i x^*_j , & \text{if } |i-j| = 1\\
    0  , & \text{otherwise }
    \end{cases}.
\end{align*}
The optimization problem \eqref{eqn:cycle} does not have any spurious solution $\mathbf{\hat{x}}$ such that $\hat{x}_i = 0$ for some $i = 0, \dots 2k$. To prove this, assume without loss of generality that $\hat{x}_0 = 0$ for a stationary point $\mathbf{\hat{x}}$. By the proof of Theorem \ref{thm: bm-solves-chain}, we know that $\hat{x}_{2k-1}, \hat{x}_{2k}, \hat{x}_{1}, \hat{x}_{2} \not = 0$. Let $m$ be defined as an index such that $\hat{x}_{1}, \hat{x}_{2}, \dots, \hat{x}_{m} \not = 0$ and either $\hat{x}_{m+1} = 0$ or $m=n$. Thus, $2 < m < 2k-2$. One can characterize the stationary points as. 
    \begin{align*}
        & [\nabla f(\mathbf{\hat{x}})]_{2k}  =  0  \xrightarrow{} \hat{x}_{2k-1} \hat{x}_{2k} = x^*_{2k-1} x^*_{2k},  \\
        & [\nabla f(\mathbf{\hat{x}})]_0  = 0 \xrightarrow{}  - x^*_0 x^*_{2k} \hat{x}_{2k} = x^*_0 x^*_1 \hat{x}_1,  \\
        & [\nabla f(\mathbf{\hat{x}})]_1  = 0  \xrightarrow{} \hat{x}_1 \hat{x}_2 = x^*_1 x^*_2,  \\
        & [\nabla f(\mathbf{\hat{x}})]_2  = 0  \xrightarrow{} \hat{x}_2 \hat{x}_3 = x^*_2 x^*_3,  \\
        & \vdots \\
        & [\nabla f(\mathbf{\hat{x}})]_m  = 0  \xrightarrow{} \hat{x}_{m-1} \hat{x}_m = x^*_{m-1} x^*_m.
    \end{align*}
    Setting $\hat{x}_{2k}$ as a free variable gives the following characterization of the stationary points:
    \begin{align*}
        & \hat{x}_{2k-2}  = \frac{x^*_{2k-2}}{x^*_{2k}}\hat{x}_{2k}, \\
        & \hat{x}_{2k-1}  = \frac{x^*_{2k-1}x^*_{2k}}{\hat{x}_{2k}}, \\
        & \hat{x}_{2k}  = \hat{x}_{2k}, \\
        & \hat{x}_0  = 0, \\
        & \hat{x}_1  = - \frac{x^*_{2k}\hat{x}_{2k}}{(x^*_1)^2} x^*_1 = -\alpha x^*_1, \\
        & \hat{x}_2  = - \frac{(x^*_1)^2}{x^*_{2k}\hat{x}_{2k}} x^*_2 = - \alpha^{-1} x^*_2, \\
        & \hat{x}_3  = -\alpha x^*_3, \\
        & \vdots \\
        & \hat{x}_m  = - \alpha^{(-1)^{m-1}} x^*_m.
    \end{align*}
    Similar to proof of Theorem \ref{thm: bm-solves-chain}, we focus on the following $(m+2) \times (m+2)$ submatrix of the Hessian that is $[\nabla^2 f(\mathbf{\hat{x}})]_{(2k-1):m, (2k-1):m}$ where $(2k-1):m$ denotes the rows/columns corresponding to $(2k-1), 2k, 0, 1, \dots, m$ in that respective order:
    \begin{align*}
        [\nabla^2 f(\mathbf{\hat{x}})]_{(2k-1):m, (2k-1):m} = \begin{pmatrix}
          \mathbf{E}
          & \rvline & \mathbf{D} \\
        \hline
          \mathbf{D}^T & \rvline &
          \begin{matrix}
          \mathbf{A} & \mathbf{B} \\
          \mathbf{B}^T & \mathbf{C}
          \end{matrix}
        \end{pmatrix},
    \end{align*}
    where $\mathbf{A} \in \mathbb{R}^{3 \times 3}, \mathbf{B} \in \mathbb{R}^{3 \times (m-1)}, \mathbf{C} \in \mathbb{R}^{(m-1) \times (m-1)}, \mathbf{D} \in \mathbb{R}^{1 \times (m+2)}$ and $\mathbf{E} \in \mathbb{R}$. The submatrices can be written as 
    \begin{align*}
        \mathbf{A} & =
        \begin{bmatrix}
            \frac{(x^*_{2k-1})^2 (x^*_{2k})^2}{(\hat{x}_{2k})^2} & -x^*_{2k} x^*_{0} & 0 \\
            - x^*_{2k} x^*_{0} & (\hat{x}_{2k})^2 + \frac{(x^*_{2k})^2 (\hat{x}_{2k})^2}{(x^*_1)^2} & - x^*_{0} x^*_{1} \\
            0 & - x^*_{0} x^*_{1} & \frac{(x^*_{1})^4 (x^*_2)^2}{(x^*_{2k})^2 (\hat{x}_{2k})^2}
        \end{bmatrix}, \\
        \mathbf{B} & = \begin{bmatrix}
            0 & 0 & \dots & 0 \\
            0 & 0 & \dots & 0 \\
            x^*_{1} x^*_{2} & 0 & \dots & 0 
        \end{bmatrix}, \\
        \mathbf{C} & = \begin{bmatrix}
            \alpha^2 (x^*_{1})^2 + \alpha^2 (x^*_{3})^2 & x^*_{2} x^*_{3} & 0 & \dots & 0 \\
            x^*_{2} x^*_{3} & \alpha^{-2} (x^*_{2})^2 + \alpha^{-2} (x^*_{4})^2 & x^*_{3} x^*_{4} & \dots & 0 \\
            0 & x^*_{3} x^*_{4} &  \alpha^2 (x^*_{3})^2 + \alpha^2 (x^*_{5})^2 & \dots & 0 \\
            \vdots & \vdots & \vdots & \ddots & \vdots \\
            0 & 0 & 0 & \dots & x^*_{m-1} x^*_{m} \\
            0 & 0 & 0 & \dots & (\alpha^2)^{({(-1)}^{m-2})} (x^*_{m-1})^2
        \end{bmatrix}, \\
        \mathbf{D} & = \begin{bmatrix}
            x^*_{2k-1} x^*_{2k} & 0 & \dots & 0
        \end{bmatrix}, \\
        \mathbf{E} & = (\hat{x}_{2k})^2 + \frac{(x^*_{2k-2})^2}{(x^*_{2k})^2} (\hat{x}_{2k})^2.
        \end{align*}
        Similar to proof of Theorem \ref{thm: bm-solves-chain}, we can investigate three different cases to demonstrate that $\mathbf{\hat{x}}$ cannot be a second-order critical point, which are $\mathbf{A} \not \succeq 0$, $\mathbf{A} \succeq 0$ with at least one zero eigenvalue, and $\mathbf{A} \succ 0$. Firstly, if $\mathbf{A} \not \succeq 0$, $\nabla^2 f(\mathbf{\hat{x}})$ cannot be positive semidefinite because $\mathbf{A}$ is a principal minor of the Hessian. The second case is when $\mathbf{A}$ is positive semidefinite but not positive definite. In that case, for a vector $\mathbf{v}$ that satisfies $\mathbf{A} \mathbf{v} = 0$, we know that $v_1$ and $v_3$ cannot be zero simultaneously. Consider following principal minors that correspond to $[\nabla^2 f(\hat{\mathbf{x}})]_{(2k-1):1, (2k-1):1}$, which includes the rows/columns corresponding to $(2k-1), 2k, 0, 1$ in that order, and $[\nabla^2 f(\hat{\mathbf{x}})]_{(2k):2, (2k):2}$, which includes the rows/columns corresponding to $2k, 0, 1, 2$ in that order, respectively:
        \begin{align*}
        \begin{pmatrix}
           (\hat{x}_{2k})^2 + \frac{(x^*_{2k-2})^2}{(x^*_{2k})^2} (\hat{x}_{2k})^2
          & \rvline & \begin{matrix}
          x^*_{2k-1} x^*_{2k} & 0 & 0
          \end{matrix}  \\
        \hline
          \begin{matrix}
           x^*_{2k-1} x^*_{2k} \\ 0 \\ 0
          \end{matrix}  & \rvline &
          \mathbf{A}
        \end{pmatrix}, \qquad  \begin{pmatrix}
           \mathbf{A}
          & \rvline & \begin{matrix}
          0 \\ 0  \\ x^*_{1} x^*_{2}
          \end{matrix}  \\
        \hline
          \begin{matrix}
          0 & 0 & x^*_{1} x^*_{2}
          \end{matrix}  & \rvline &
          \alpha^2 (x^*_{1})^2 + \alpha^2 (x^*_{3})^2
        \end{pmatrix}.
    \end{align*}
     The submatrices are not positive semidefinite by Schur complement if
    \begin{align*}
        \mathbf{A} - \frac{1}{(\hat{x}_{2k})^2 + \frac{(x^*_{2k-2})^2}{(x^*_{2k})^2} (\hat{x}_{2k})^2} \begin{bmatrix}
            (x^*_{2k-1})^2 (x^*_{2k})^2 & 0 & 0 \\ 0 & 0 & 0 \\ 0 & 0 & 0
        \end{bmatrix} \not \succeq 0
    \end{align*}
    and
    \begin{align*}
        \mathbf{A} - \frac{1}{\alpha^2 (x^*_{1})^2 + \alpha^2 (x^*_{3})^2} \begin{bmatrix}
            0 & 0 & 0 \\ 0 & 0 & 0 \\ 0 & 0 & (x^*_{1})^2 (x^*_{2})^2
        \end{bmatrix} \not \succeq 0.
    \end{align*}
    If $v_1 \not = 0$, then 
    \begin{align*}
        \mathbf{v}^T \mathbf{A} \mathbf{v} - \mathbf{v}^T \frac{1}{(\hat{x}_{2k})^2 + \frac{(x^*_{2k-2})^2}{(x^*_{2k})^2} (\hat{x}_{2k})^2} \begin{bmatrix}
            (x^*_{2k-1})^2 (x^*_{2k})^2 & 0 & 0 \\ 0 & 0 & 0 \\ 0 & 0 & 0
        \end{bmatrix} \mathbf{v} = 0 - \frac{v_1^2 (x^*_{2k-1})^2 (x^*_{2k})^2}{(\hat{x}_{2k})^2 + \frac{(x^*_{2k-2})^2}{(x^*_{2k})^2} (\hat{x}_{2k})^2} < 0.
    \end{align*}
    Hence, the quadratic form of the Hessian has a descent direction. If $v_1 = 0$, then $v_3 \not =0$ and we can use the second submatrix to show that the submatrix is not positive semidefinite. Thus, $\mathbf{\hat{x}}$ is not a second-order critical point. As a result, if $\mathbf{\hat{x}}$ is a second-order critical point, then $\mathbf{A}$ must be positive definite. In that case, however, we can show that the submatrix 
    \begin{align*}
        \begin{bmatrix}
            \mathbf{A} & \mathbf{B} \\
            \mathbf{B}^T & \mathbf{C}
        \end{bmatrix} 
    \end{align*}
    cannot be positive semidefinite, which implies that $\mathbf{\hat{x}}$ cannot be a second-order critical point.  We prove the last proposition by showing that $\mathbf{C} - \mathbf{B}^T \mathbf{A}^{-1} \mathbf{B} \not \succeq 0$ by using Schur complement idea. Note that 
    \begin{align*}
        \mathbf{B}^T \mathbf{A}^{-1} \mathbf{B} = \begin{bmatrix}
            \mathbf{A}^{-1}_{3,3} (x^*_{1})^2 (x^*_{2})^2 & 0 & \dots & 0 \\
            0 & 0 & \dots & 0 \\
            \vdots & \vdots & \ddots & \vdots \\
            0 & 0 & \dots & 0
        \end{bmatrix}.
    \end{align*}
    We calculate $\mathbf{A}_{3,3}^{-1}$ by using the cofactors of $\mathbf{A}$ as follows: 
    \begin{align*}
        \mathbf{A}_{3,3}^{-1} = \frac{1}{\text{det}(\mathbf{A})} \Big( (x^*_{2k-1})^2 (x^*_{2k})^2 + \frac{(x^*_{2k-1})^2 (x^*_{2k})^4}{(x^*_{1})^2} - (x^*_{2k})^2 (x^*_{0})^2 \Big),
    \end{align*}
    where 
    \begin{align*}
        \text{det}(\mathbf{A}) = \frac{(x^*_{2k-1})^2 (x^*_{2k})^2}{(\hat{x}_{2k})^2} \Big( \frac{(x^*_{1})^4 (x^*_{2})^2}{(x^*_{2k})^2} + (x^*_{1})^2 (x^*_{2})^2 - (x^*_{0})^2 (x^*_{1})^2 \Big) - \frac{(x^*_{0})^2 (x^*_{1})^4 (x^*_{2})^2}{(\hat{x}_{2k})^2}.
    \end{align*}
    An algebraic manipulation shows that $\frac{d \mathbf{A}_{3,3}^{-1}}{d (x^*_0)^2} > 0$. Thus, the value of $\mathbf{A}_{3,3}^{-1}$ is minimized when $(x^*_0)^2 \xrightarrow{} 0$. Whenever $x^*_0 = 0$, $\mathbf{A}^{-1}_{3,3} (x^*_{1})^2 (x^*_{2})^2$ is equal to $\alpha^2 (x^*_{1})^2$. As a result, $\mathbf{A}^{-1}_{3,3} (x^*_{1})^2 (x^*_{2})^2 > \alpha^2 (x^*_{1})^2$ by the non-zero assumption on $u_0$. We can write the matrix $\mathbf{C} - \mathbf{B}^T \mathbf{A}^{-1} \mathbf{B}$ as the summation of $m-1$ matrices as
    \begin{align*}
        \mathbf{C} - \mathbf{B}^T \mathbf{A}^{-1} \mathbf{B} = & \begin{bmatrix}
            \alpha^2 (x^*_{1})^2 - \mathbf{A}_{3,3}^{-1} (x^*_{1})^2 (x^*_{2})^2 & 0 & \dots & 0 & 0 \\
            0 & 0 & \dots & 0 & 0 \\
            \vdots & \vdots & \ddots & \vdots & \vdots  \\
            0 & 0 & \dots & 0 & 0 \\
            0 & 0 & \dots & 0 & 0
        \end{bmatrix} + \\
        & \begin{bmatrix}
            \alpha^2 (x^*_{3})^2 & x^*_{2} x^*_{3} & \dots & 0 & 0 \\
            x^*_{2} x^*_{3} & \alpha^{-2} (x^*_{2})^2 & \dots & 0 & 0 \\
            \vdots & \vdots & \ddots & \vdots & \vdots  \\
            0 & 0 & \dots & 0 & 0 \\
            0 & 0 & \dots & 0 & 0
        \end{bmatrix} + \cdots + \\
        & \begin{bmatrix}
            0 & 0 & \dots & 0 & 0 \\
            0 & 0 & \dots & 0 & 0 \\
            \vdots & \vdots & \ddots & \vdots & \vdots  \\
            0 & 0 & \dots & (\alpha^2)^{(-1)^{m-1}} (x^*_m)^2 & x^*_{m-1} x^*_m \\
            0 & 0 & \dots & x^*_{m-1} x^*_m & (\alpha^2)^{(-1)^{m-2}} (x^*_{m-1})^2 
        \end{bmatrix}.
    \end{align*}
    Consider the vector $\mathbf{v}$ defined as $v_{2t+1} = \frac{x^*_{2t + 2}}{x^*_{2}}$ and $v_{2t} = \alpha^2 \frac{x^*_{2t+1}}{x^*_{2}}$ for $t = 0, \dots, \lfloor (m-1)/2 \rfloor$ so that $\mathbf{v}^T (\mathbf{C} - \mathbf{B}^T \mathbf{A}^{-1} \mathbf{B}) \mathbf{v}$ results in zero for all the matrices above except for the first one. Then,
    \begin{align*}
        \mathbf{v}^T (\mathbf{C} - \mathbf{B}^T \mathbf{A}^{-1} \mathbf{B}) \mathbf{v} = \alpha^2 (x^*_{1})^2 - \mathbf{A}_{3,3}^{-1} (x^*_{1})^2 (x^*_{2})^2 + 0 + \dots + 0 < 0.
    \end{align*}
    As a result, $ \mathbf{C} - \mathbf{B}^T \mathbf{A}^{-1} \mathbf{B}$ is not positive definite. Hence, the Hessian is not positive definite either. Consequently, the problem cannot have a second-order critical point $\mathbf{\hat{x}}$ with $\hat{x}_i^* = 0$ for some $i \in [n]$.
    
    Any solution that meets the requirements of second-order necessary optimality conditions cannot contain any zero entries. In addition, if a stationary point $\mathbf{\hat{x}}$ satisfies $\hat{x}_i \hat{x}_j = x^*_i x^*_j$ for some $(i,j) \in \mathcal{E}$, then $\hat{x}_i \hat{x}_j = x^*_i x^*_j$ for all $(i,j) \in \mathcal{E}$ due to the stationarity condition. Therefore, a spurious solution has the following properties: $\hat{x}_i \hat{x}_j \not = x^*_i x^*_j$ for all $(i,j) \in \mathcal{E}$ and $\hat{x}_i \not = 0$ for all $i = 0, \dots 2k$. 

Define $a_{i,j} = x_i x_j - x^*_i x^*_j$. Then, the stationarity condition for the B-M factorized problem \eqref{eqn:cycle} can be written as 
\begin{equation*}
    [\nabla f(\mathbf{\hat{x}})]_i = \begin{cases}
    \hat{a}_{0,2k} \hat{x}_{2k} + \hat{a}_{0,1} \hat{x}_1 = 0 , & \text{if }i =0\\
    \hat{a}_{i-1,i} \hat{x}_{i-1} + \hat{a}_{i,i+1} \hat{x}_{i+1} = 0 , & \text{if }i \not \in  \{0,2k\} \\
    \hat{a}_{2k-1,2k} \hat{x}_{2k-1} + \hat{a}_{0,2k}\hat{x}_0 = 0, & \text{if }i =2k
    \end{cases}.
\end{equation*}
The following calculation yields a contradiction because none of the terms can be zero:
\begin{equation*}
    \sum_{i=0}^{2k} (-1)^i [\nabla f(\mathbf{\hat{x}})]_i = 2 \hat{a}_{0,2k} \hat{x}_0 \hat{x}_{2k} = 0.
\end{equation*}
As a result, all the second-order critical points are global solutions and we obtain the valid factors of the ground truth matrix completion. 
\end{proof}

\subsection{Proof of Theorem \ref{thm: sdp-fails-cycle}}

\begin{proof}
     We aim to find an instance of the matrix completion problem with an odd-numbered cycle graph for which the SDP problem \eqref{eqn:ms-mc-sdp} fails. Consider the following rank-1 feasible solution $\mathbf{\hat{M}} = \mathbf{yy}^T$ with
\[ y_0 = (\mathbf{\hat{M}}_{0,0})^{1/2} \geq x^*_0, \quad y_{2t-1} = \left(\frac{(x^*_{2t-1})^2 (x^*_0)^2}{\mathbf{\hat{M}}_{0,0}}\right)^{1/2}, \quad y_{2t} = \left(\frac{(x^*_{2t})^2}{(x^*_0)^2}\mathbf{\hat{M}}_{0,0}\right)^{1/2}, \]
where $t = 1, \dots, k$. Note that $\mathbf{\hat{M}}_{0,0}$ is free to choose and we can minimize the trace to find the best rank-1 solution. This is equivalent to the following optimization problem:
\begin{equation*}
    \min_{\mathbf{\hat{M}}_{0,0} \geq (x^*_0)^2} \frac{\mathbf{\hat{M}}_{0,0}}{(x^*_0)^2} \sum_{t=0}^{k} (x^*_{2t})^2 + \frac{(x^*_0)^2}{\mathbf{\hat{M}}_{0,0}} \sum_{t=1}^{k} (x^*_{2t-1})^2.
\end{equation*}
Note that whenever $\mathbf{\hat{M}}_{0,0} = (x^*_0)^2$, the solution $\mathbf{\hat{M}}$ is the same as ground truth solution. A basic first-order optimality condition implies that whenever $\sum_{t=1}^{k} (x^*_{2t-1})^2 > \sum_{t=0}^{k} (x^*_{2t})^2$, the optimal $\mathbf{\hat{M}}_{0,0}$ is strictly larger than $ (x^*_0)^2$. However, the objective value is strictly better than the trace of the ground truth matrix. Due to the symmetry of the problem, any node can be chosen as node 0. Therefore, if the condition $\sum_{t=1}^{k} (x^*_{2t-1})^2 > \sum_{t=0}^{k} (x^*_{2t})^2$ holds for some chosen node $0$, the SDP fails to solve the matrix completion problem.
\end{proof}

\end{document}